\pdfoutput=1
\RequirePackage{ifpdf}
\ifpdf % We are running pdfTeX in pdf mode
\documentclass[pdftex]{sigma}
\else
\documentclass{sigma}
\fi

\begin{document}

\allowdisplaybreaks

\renewcommand{\thefootnote}{$\star$}

\renewcommand{\PaperNumber}{030}

\FirstPageHeading

\ShortArticleName{Motions of Curves in the Projective Plane Inducing the Kaup--Kupershmidt Hierarchy}

\ArticleName{Motions of Curves in the Projective Plane \\ Inducing the Kaup--Kupershmidt Hierarchy\footnote{This
paper is a contribution to the Special Issue ``Symmetries of Dif\/ferential Equations: Frames, Invariants and Applications''. The full collection is available at \href{http://www.emis.de/journals/SIGMA/SDE2012.html}{http://www.emis.de/journals/SIGMA/SDE2012.html}}}

\Author{Emilio~MUSSO}

\AuthorNameForHeading{E.~Musso}

\Address{Dipartimento di Scienze Matematiche, Politecnico di Torino,\\
Corso Duca degli Abruzzi 24, I-10129 Torino, Italy}
\Email{\href{mailto:emilio.musso@polito.it}{emilio.musso@polito.it}}

\ArticleDates{Received February 08, 2012, in f\/inal form May 11, 2012; Published online May 24, 2012}

\Abstract{The equation of a motion of curves in the projective plane is deduced. Local f\/lows are def\/ined in terms of polynomial dif\/ferential functions. A family of local f\/lows inducing the Kaup--Kupershmidt hierarchy is constructed. The integration of the congruence curves is discussed. Local motions def\/ined by the traveling wave cnoidal solutions of the f\/ifth-order Kaup--Kupershmidt equation are described.}


\Classification{53A20; 53A55; 33E05; 35Q53; 37K10}

%%%%%%%%%%%%%%%%%%%%%
\def\C{\mathbb{C}}
\def\R{\mathbb{R}}
\def\RP{\mathbb{RP}}
\def\Z{\mathbb{Z}}
\def\E{\mathbb{E}}
%%%%%%%%%%%%%%%%%%%%%%%%%%%%

\renewcommand{\thefootnote}{\arabic{footnote}}
\setcounter{footnote}{0}

\section{Introduction}\label{s:intro}
The interrelations between hierarchies of integrable non linear evolution equations and motions of curves have been widely investigated in the last decades, both in geometry and mathematical physics. In the seminal papers \cite{GP1,GP,NSW}, Goldstein, Petrich and Nakayama, Segur, Wadati, showed that the mKdV hierarchy can be deduced from local motions of curves in the Euclidean plane. Later, this result was extended to other $2$-dimensional geometries \cite{CQ3,CQ1,CQ2,CQ5,CQ4,Pi,QSL} or to higher-dimensional homogeneous spaces \cite{AB,CIB-PhD,HS-PAMS02,Ivey-cpm,LP,LQS,MN2}. The invariant curve f\/lows related to integrable hierarchies are induced by inf\/inite-dimensional Hamiltonian systems def\/ined by invariant functionals and geometric Poisson brackets on the space of dif\/ferential invariants of parameterized curves \cite{MB1,MB2}. Another feature is the existence of f\/inite-dimensional reductions leading to Liouville-integrable Hamiltonian systems. Typically, these reductions correspond to curves which evolve by congruences of the ambient space. They have both a variational and a Hamiltonian description: as extremals of a geometric variational problem def\/ined by the conserved densities of the hierarchy, and as solutions of a f\/inite-dimensional integrable contact Hamiltonian systems. Examples include local motions of curves in two-dimensional Riemannian space-forms \cite{M1}, local motions of star-shaped curves in centro-af\/f\/ine geometry \cite{M2,Pi} and local motions of null curves in $3$-dimensional pseudo-Riemannian space forms \cite{MN2}. In \cite{CQ3, CQ1,CQ2,CQ5}, K.-S.~Chu and C.~Qu gave a rather complete account of the integrable hierarchies originated by local motions of curves in $2$-dimensional Klein geometries. In particular, they showed that the f\/ifth-order Kaup--Kupershmidt equation is induced by a local motion in centro-af\/f\/ine geometry and that modif\/ied versions of the f\/ifth- and seventh-order Kaup--Kupershmidt equations can be related to local motions of curves in the projective plane. Our goal is to demonstrate the existence of a sequence of local motions of curves in projective plane inducing the entire Kaup--Kupershmidt hierarchy. We analyze the congruence curves of the f\/lows and we investigate in more details the congruence curves associated to the cnoidal traveling wave solutions \cite{Za} of the f\/ifth-order Kaup--Kupershmidt equation. In particular, we show that the critical curves of the projective invariant functional \cite{Ca,MN1} are congruence curves of the f\/irst f\/low of the hierarchy.

The material is organized into three sections. In the f\/irst section we collect basic facts about the Kaup--Kupershmidt hierarchy from the existing literature \cite{FG,FuOe,Ku,MV1,CR,MV2,We} and we discuss an alternative form of the equations of the hierarchy. We examine the cnoidal traveling wave solutions of the f\/ifth-order Kaup--Kupershmidt equation \cite{Za} and we exhibit new traveling wave solutions in terms of Weierstrass elliptic functions. In the second section we recall the construction of the projective frame along a plane curve without inf\/lection or sextatic points and we introduce the projective line element and the projective curvature \cite{Ca,Ha,OT,Wi}. Subsequently, we use the Cartan's canonical frame \cite{Ca,MN1} to study the equations of a motion of curves in the projective plane. Consequently, we construct local motions in terms of dif\/ferential polynomial functions and we deduce the existence of a sequence of local f\/lows inducing the Kaup--Kupershmidt hierarchy. In the third section we focus on congruence motions of the f\/lows. In particular, we explicitly implement the general process of integration to determine the motions of the critical curves of the projective invariant functional \cite{Ca,MN1}. The symbolical and numerical computations as well the graphics have been worked out with the software {\it Mathematica 8}. We adopt \cite{L} as a reference for elliptic functions and integrals.

\section[The Kaup-Kupershmidt hierarchy]{The Kaup--Kupershmidt hierarchy}\label{s:1}

\subsection{Preliminaries and notations}\label{ss:1.1}

Let $J(\R,\R)$ be the jet space of smooth functions $u:\R\to \R$, equipped with the usual coordinates
$(s,u_{(0)}, u_{(1)},\dots, u_{(h)},\dots)$. The prolongation of a smooth function
$u$ is denoted by  $j(u)$. Similarly, if $u(s,t)$ is a function of the variables $s$ and $t$, its partial prolongation with respect to the $s$-variable will be denoted by $j_s(u)$. A map $\mathfrak p : J(\R,\R) \to \R$ is said to be a \textit{polynomial differential
function} if
\[
\mathfrak{p}(\mathbf{u})=P(u_{(0)},u_{(1)},\dots,u_{(h)})\qquad \forall\, \mathbf{u} \in J(\R,\R),
\]
where $P$ is a polynomial in $h+1$ variables.
The algebra of polynomial dif\/ferential functions, $J[\mathbf{u}]$,
is endowed with the {\it total derivative}
\[
D \mathfrak p =
\sum_{p=0}^{\infty} \frac{\partial p}{\partial u_{(p)}} u_{(p+1)}
\]
and the {\it Euler operator}
\[
 \mathcal{E}(\mathfrak p) =\sum _{\ell =0}^\infty (-1)^\ell D^\ell
 \left(\frac{\partial \mathfrak p}{\partial u_{(\ell)}}\right).
\]
For each $\mathfrak{p}\in \mathrm{Ker}(\mathcal{E})$ there is a unique $D^{-1} \mathfrak p\in J[\mathbf{u]}$ such that
\[
\mathfrak p=D\big( D^{-1}\left(\mathfrak p\right)\big),\qquad D^{-1} \mathfrak p|_0=0.
\]
We consider the linear subspace
\begin{gather}
P[\mathbf{u}]=\big\{\mathfrak{p}\in J[\mathbf{u}] : \; \mathcal{E}\big(u_{(0)}D^3\mathfrak{p}+4u_{(0)}^2D\mathfrak{p}\big)=0\big\}.\label{P[u]}
\end{gather}
The image of $P[\mathbf{u}]$ by the total derivative is denoted by  $P'[\mathbf{u}]\subset J[\mathbf{u}]$. Next we consider the integro-dif\/ferential operators
\begin{gather*}
\Delta(\phi,u) =\phi_{3s}+2u\phi_s+u_s\phi,\\
\Xi(\phi,u) =\phi_{3s}+8u\phi_s+7u_s\phi+2\big(u_{2s}+4u^2\big)\int_0^s \phi dr+2\int_0^s\big(u\phi_{2s}+4u^2\phi\big)dr,\\
\Theta(\phi,u)  = \frac{u_s}{9}\int_0^s \big(u\phi_{3s}+4u^2\phi_s\big)dr+\frac{1}{9}\big(20u^2u_s+25u_su_{2s}+10uu_{3s}+u_{5s}\big)\phi\\
\hphantom{\Theta(\phi,u)  =}{}  +\left(\frac{3}{2}+\frac{16}{9}u^3+\frac{71}{18}u_s^2+\frac{41}{9}uu_{2s}+\frac{13}{18}u_{4s}\right)\phi_s+
\left(\frac{59}{9}uu_s+\frac{35}{18}u_{3s}\right)\phi_{2s}\\
\hphantom{\Theta(\phi,u)  =}{}
 +\left(2u^2+\frac{49}{18}u_{2s}\right)\phi_{3s}-2u_s\phi_{4s}+\frac{2}{3}u\phi_{5s}+\frac{1}{18}\phi_{7s}
\end{gather*}
and the linear operators
\[
\mathcal{D}: \ J[\mathbf{u}]\to J[\mathbf{u}],\qquad \mathcal{J}:\ P'[\mathbf{u}]\to J[\mathbf{u}],\qquad \mathcal{S}: \ P[\mathbf{u}]\to J[\mathbf{u}]
\]
def\/ined by
 \begin{gather*}
 \mathcal{D}(\mathfrak{w})   = D^3\mathfrak{w} + 2u_{(0)}D\mathfrak{w} + u_{(1)}\mathfrak{w},\\
 \mathcal{J}(\mathfrak{q})  = D^{3}\mathfrak{q}+8u_{(0)}D\mathfrak{q}+7u_{(1)}\mathfrak{q}+2\big(u_{(1)}+4u_{(0)}^2\big)D^{-1}\mathfrak{q}+
 2D^{-1}\big(u_{(0)}D^{2}\mathfrak{q}+4u_{(0)}^2\mathfrak{q}\big),\\
 \mathcal{S}(\mathfrak{p}) =\frac{1}{9}u_{(1)} D^{-1}\big(u_{(0)}D^3\mathfrak{p}+4u_{(0)}^2D\mathfrak{p}\big)+\frac{1}{9}\big(20u_{(0)}^2u_{(1)}+25u_{(1)}u_{(2)}+10u_{(0)}u_{(3)}+u_{(5)}\big)\mathfrak{p}\\
\hphantom{\mathcal{S}(\mathfrak{p}) =}{}
 +\left(\frac{3}{2}+\frac{16}{9}u_{(0)}^3+\frac{71}{18}u_{(1)}^2+\frac{41}{9}u_{(0)}u_{(2)}+\frac{13}{18}u_{(4)}\right)D\mathfrak{p}+
\left(\frac{59}{9}u_{(0)}u_{(1)}+\frac{35}{18}u_{(3)}\right)D^2\mathfrak{p}\\
\hphantom{\mathcal{S}(\mathfrak{p}) =}{}
 +\left(2u_{(0)}^2+\frac{49}{18}u_{(2)}\right)D^3\mathfrak{p}-2u_{(1)}D^4\mathfrak{p}+\frac{2}{3}u_{(0)}D^5\mathfrak{p}+\frac{1}{18}D^7\mathfrak{p}.
\end{gather*}
From the def\/inition of the operators is clear that
\[
\mathcal{D}(\mathfrak{w})|_{j(u)}=\Delta(\mathfrak{w}|_{j(u)},u),\qquad
\mathcal{S}(\mathfrak{p})|_{j(u)}=\Theta(\mathfrak{p}|_{j(u)},u),\qquad \mathcal{J}(\mathfrak{q})|_{j(u)}=\Xi(\mathfrak{q}|_{j(u)},u),
\]
for every $\mathfrak{w}\in j[\mathbf{u}]$, $\mathfrak{p}\in P[\mathbf{u}]$, $\mathfrak{q}\in P'[\mathbf{u}]$ and every $u\in C^{\infty}(\R,\R)$.

\begin{lemma}\label{lemma1}
The operators $\mathcal{D}$, $\mathcal{J}$ and $\mathcal{S}$ satisfy
\begin{gather*}%\label{1:1.11}
\mathcal{D}\mathcal{J}(\mathfrak{p})=18\mathcal{S}\big(D^{-1}\mathfrak{p}\big)-27\mathfrak{p}\qquad \forall\, \mathfrak{p}\in P'[\mathbf{u}].
\end{gather*}
\end{lemma}

\begin{proof} A direct computation shows that
\[\Delta(\Xi(\phi,u),u)= 18\Theta\left(\int_0^s\phi dr,u\right)-27\phi \qquad \forall\,  u,\phi\in C^{\infty}(\R,\R).
\]
This implies
 \begin{gather*}
 \mathcal{D}\mathcal{J}(\mathfrak{p})|_{j(u)} =\Delta(\Xi(\mathfrak{p}|_{j(u)},u),u)=
18\Theta(D^{-1}(\mathfrak{p})|_{j(u)},u)-27\mathfrak{p}|_{j(u)}\\
\hphantom{\mathcal{D}\mathcal{J}(\mathfrak{p})|_{j(u)}}{}
 = 18\mathcal{S}(D^{-1}(\mathfrak{p}))|_{j(u)}-27\mathfrak{p}|_{j(u)},
 \end{gather*}
for every $\mathfrak{p}\in P[\mathbf{u}]$ and every $u\in C^{\infty}(\R,\R)$. This yields the required result.
 \end{proof}

\subsection{Construction of the hierarchy}\label{ss:1.2}

According to \cite{FuOe,CR,We} there are three sequences
\[\{\mathfrak h_n\}_{n\in \mathbf{N}}\subset \mathrm{Im}(\mathcal{E}),\qquad \{\mathfrak q_n\}_{n\in \mathbf{N}}\subset P[\mathbf{u}],\qquad \{\mathfrak p_n\}_{n\in \mathbf{N}}\subset J[\mathbf{u}]\]
of polynomial dif\/ferential functions
def\/ined by the recursion formulae
\begin{gather}\label{1:2.1}
\mathfrak h_{n+2}= \mathcal{J} \mathcal{D}(\mathfrak h_{n}),\qquad   \mathcal{D}(\mathfrak h_{n})=D(\mathfrak q_{n}),\qquad
\mathfrak h_{n}=\mathcal{E}(\mathfrak{p}_n)
\end{gather}
and by the initial data
\begin{gather*}
  \mathfrak h_0=1, \qquad \mathfrak h_1=u_{(2)}+4u_{(0)}^2.
\end{gather*}

\begin{definition}
The
{\it Kaup--Kupershmidt hierarchy} $\{\mathcal{K}_n\}$ is the sequence of evolution equations def\/ined by
\begin{gather*}%\label{1:2.2}
\mathcal{K}_n : \ u_t + \mathcal{D} \mathfrak h_{n}|_{j(u)} =0.
\end{gather*}
\end{definition}

  In view of (\ref{1:2.1}), the equations of the hierarchy can be written either in the Hamiltonian form
\begin{gather*}%\label{1:2.3}
\mathcal{K}_n : \ u_t+\mathcal{D}\mathcal{E}(\mathfrak p_{n})|_{j(u)}=0,
\end{gather*}
or else in the conservation form
\begin{gather*}%\label{1:2.4}
\mathcal{K}_n : \ u_t+D\mathfrak{q}_n=0.
\end{gather*}

\begin{remark}
The polynomial dif\/ferential functions $\mathfrak h_n$, $\mathfrak q_n$, $\mathfrak p_n$ and the equations of the hierarchy can be computed with any software of symbolic calculus (see Appendix~\ref{a1}). For $n=1,2$ we f\/ind
 \begin{gather*}
\mathfrak{h}_1(\mathbf{u}) =u_{(2)}+4 u_{(0)}^2,\\
\mathfrak{h}_2(\mathbf{u}) =12 u_{(2)} u_{(0)}+6 u_{(1)}^2+u_{(4)}+\frac{32}{3} u_{(0)}^3,\\
\mathfrak{q}_1(\mathbf{u}) =\frac{20}{3} u_{(0)}^3+\frac{15}{2} u_{(1)}^2+10 u_{(0)} u_{(2)}+u_{(4)},\\
\mathfrak{q}_2(\mathbf{u}) =\frac{56}{3} u_{(0)}^4+70 u_{(0)} u_{(1)}^2+56 u_{(0)}^2 u_{(2)}+\frac{49}{2} u_{(2)}^2+35 u_{(1)} u_{(3)}+14
u_{(0)} u_{(4)}+u_{(6)},\\
\mathfrak{p}_1(\mathbf{u}) =\frac{1}{2}u_{(0)}2u_{(2)}+\frac{4}{3}u_{(0)}^2,\\
\mathfrak{p}_2(\mathbf{u}) =\frac{1}{2}u_{(0)}u_{(4)}+4u_{(0)}^2u_{(2)}+2u_{(0)}u_{(1)}^2+\frac{8}{3}u_{(0)}^4.
\end{gather*}
 Consequently, the f\/irst two equations of the hierarchy are
 \begin{gather*}
 u_{t}+10uu_{3s}+25u_su_{2s}+20u^2u_s + u_{5s} =0,\\
 u_{t} + u_{7s}+14uu_{5s}+49u_su_{4s}+84u_{2s}u_{3s}+252uu_su_{2s}+70u_s^3
+56u^2u_{3s}+\frac{224}{3}u^3u_s =0.
\end{gather*}
\end{remark}

\begin{definition}
Denoting by $[r]$ the integer part of $r$, we set
\begin{gather*}%\label{1:2.5}
\ell_n = \left[\frac{n}{2}\right]-\frac{1}{2}(1+(-1)^n),\qquad \lambda_n=\frac{1}{2}(1+(-1)^n))(-27)^{\left[\frac{n}{2}\right]}
\end{gather*}
and we def\/ine $\{\mathfrak{v}_n\}_{n\in \mathbb{N}}\subset P[\mathbf{u}]$ by
\begin{gather*}%\label{1:2.6}
\mathfrak{v}_n=18\sum_{h=0}^{\ell_n}(-27)^h\mathfrak{w}_{n-2h},\qquad n = 1,\dots,
\end{gather*}
where $\mathfrak{w}_0=0$, $\mathfrak{w}_1=1/2$ and $\mathfrak{w}_n=\mathfrak{q}_{n-2}$, $n>1$.
\end{definition}

\begin{proposition}\label{proposition1}
The equations of the hierarchy can be written in the form
\begin{gather}\label{1:2.7}
\mathcal{K}_n : \ \partial_tu + \mathcal{S}(\mathfrak{v}_{n})|_{j_s(u)}+\lambda_n u_{s}=0.\end{gather}
\end{proposition}

\begin{proof}
For $n=1,2$ the proposition can be checked by a direct computation.
We prove (\ref{1:2.7}) when $n$ is odd. By induction, suppose that (\ref{1:2.7}) is true for $n=2p-1$. Note that
\[\lambda_{2p-1}=\lambda_{2p+1}=0,\qquad \ell_{2p+1}=\ell_{2p-1}+1.\]
By the inductive hypothesis we have
\[\mathcal{K}_{2p-1}(u)=u_t+D(\mathfrak{q}_{2p-1})|_{j_s(u)}=u_t+\mathcal{S}(\mathfrak{v}_{2p-1})|_{j_s(u)}=
u_t+\Theta(\mathfrak{v}_{2p-1}|_{j_s(u)},u),\]
which implies
\[D(\mathfrak{q}_{2p-1})|_{j_s(u)}=\Theta(\mathfrak{v}_{2p-1}|_{j_s(u)},u).\]
Using Lemma \ref{lemma1} we f\/ind
\begin{gather*}
\mathcal{K}_{2p+1}(u) =u_t+\mathcal{D}\mathcal{J}(D\mathfrak{q}_{2p-1})|_{j_s(u)}=
u_t+18\mathcal{S}(\mathfrak{q}_{2p-1})|_{j_s(u)}-27D\mathfrak{q}_{2p-1}|_{j_s(u)} \\
\hphantom{\mathcal{K}_{2p+1}(u)}{}
 =u_t+18\Theta(\mathfrak{q}_{2p-1}|_{j_s(u)},u)-27\Theta(\mathfrak{v}_{2p-1}|_{j_s(u)},u) \\
 \hphantom{\mathcal{K}_{2p+1}(u)}{}
 =u_t+\Theta((18\mathfrak{w}_{2p+1}-27\mathfrak{v}_{2p-1})|_{j_s(u)},u).
\end{gather*}
Using
 \begin{gather*}
18\mathfrak{w}_{2p+1}-27\mathfrak{v}_{2p-1} = 18\left(\mathfrak{w}_{2p+1}-27\sum_{h=0}^{\ell_{2p-1}}(-27)^h\mathfrak{w}_{2p-1-2h}\right) \\
\hphantom{18\mathfrak{w}_{2p+1}-27\mathfrak{v}_{2p-1}}{}
 =18\left(\mathfrak{w}_{2p+1}-\sum_{h=0}^{\ell_{2p+1}-1}(-27)^{h+1}\mathfrak{w}_{2p-1-2(h+1)}\right) \\
\hphantom{18\mathfrak{w}_{2p+1}-27\mathfrak{v}_{2p-1}}{}
 =18\left(\mathfrak{w}_{2p+1}-\sum_{h=1}^{\ell_{2p+1}}(-27)^{h}\mathfrak{w}_{2p+1-2h}\right) \\
\hphantom{18\mathfrak{w}_{2p+1}-27\mathfrak{v}_{2p-1}}{}
 =18\sum_{h=0}^{\ell_{2p+1}}(-27)^{h}\mathfrak{w}_{2p+1-2h}=\mathfrak{v}_{2p+1}
\end{gather*}
we obtain
\[
\mathcal{K}_{2p+1}(u)=u_t+\Theta(\mathfrak{v}_{2p+1}|_{j_s(u)},u)
= u_t+\mathcal{S}(\mathfrak{v}_{2p+1})|_{j_s(u)}+\lambda_{2p+1}u_s.
\]
Next we prove (\ref{1:2.7}) when $n$ is even. By induction, suppose that (\ref{1:2.7}) is true for $n=2p$.  Note that
\[
\lambda_{2p+2}=(-27)^{p+1} =-27\lambda_{2p},\qquad \ell_{2p+2}=p=\ell_{2p}+1.
\]
By the inductive hypothesis we have
 \begin{gather*}
 \mathcal{K}_{2p}(u) =u_t+D(\mathfrak{q}_{2p})|_{j_s(u)}=u_t+\mathcal{S}(\mathfrak{v}_{2p})|_{j_s(u)}-\lambda_{2p}u_s
 = u_t+\Theta(\mathfrak{v}_{2p}|_{j_s(u)},u)+\lambda_{2p}u_s,
\end{gather*}
which implies
\[
D(\mathfrak{q}_{2p})|_{j_s(u)}=\Theta(\mathfrak{v}_{2p}|_{j_s(u)},u)+\lambda_{2p}u_s.
\]
From  Lemma \ref{lemma1} we have
\begin{gather*}
\mathcal{K}_{2p+2}(u) =u_t+\mathcal{D}\mathcal{J}(D\mathfrak{q}_{2p})|_{j_s(u)}=
u_t+18\mathcal{S}(\mathfrak{q}_{2p})|_{j_s(u)}-27D\mathfrak{q}_{2p}|_{j_s(u)} \\
\phantom{\mathcal{K}_{2p+2}(u)}{}
=u_t+18\Theta(\mathfrak{q}_{2p}|_{j_s(u)},u)-27(\Theta(\mathfrak{v}_{2p}|_{j_s(u)},u)+\lambda_{2p}u_s) \\
\phantom{\mathcal{K}_{2p+2}(u)}{}
=u_t+\Theta((18\mathfrak{w}_{2p+2}-27\mathfrak{v}_{2p})|_{j_s(u)},u)+\lambda_{2p+2}u_s.
\end{gather*}
Using
 \begin{gather*}
18\mathfrak{w}_{2p+2}-27\mathfrak{v}_{2p} = 18\left(\mathfrak{w}_{2p+2}-27\sum_{h=0}^{\ell_{2p}}(-27)^h\mathfrak{w}_{2p-2h}\right) \\
\qquad{}  =18\left(\mathfrak{w}_{2p+2}-\sum_{h=0}^{\ell_{2p+2}-1}(-27)^{h+1}\mathfrak{w}_{2p-2h}\right)
 =18\left(\mathfrak{w}_{2p+2}-\sum_{h=1}^{\ell_{2p+2}}(-27)^{h}\mathfrak{w}_{2p+2-2h}\right) \\
 \qquad{}  =18\sum_{h=0}^{\ell_{2p+2}}(-27)^{h}\mathfrak{w}_{2p+2-2h}=\mathfrak{v}_{2p+2}
\end{gather*}
we f\/ind
\begin{gather*}
\mathcal{K}_{2p+2}(u)=u_t+\Theta(\mathfrak{v}_{2p+1}|_{j_s(u)},u))
+\lambda_{2p+2}u_s= u_t+\mathcal{S}(\mathfrak{v}_{2p+2})|_{j_s(u)}+\lambda_{2p+2}u_s.\tag*{\qed}
\end{gather*}
\renewcommand{\qed}{}
\end{proof}

\subsection[Traveling waves of the fifth-order Kaup-Kupershmidt equation]{Traveling waves of the f\/ifth-order Kaup--Kupershmidt equation}\label{ss:1.3}

Several classes of traveling wave solutions of the f\/ifth-order equation $\mathcal{K}_1$ have been considered in the literature \cite{Ka,Wa,Za}. In this section we generalize the elliptic families examined in \cite{Za}. The hyperbolic traveling waves found in \cite{Wa} can be obtained as limiting cases, when the parameter of the elliptic functions tends to $1$. First consider the third-order ODE
\begin{gather}\label{1:3.1}
k'''+8k  k'=0,
\end{gather}
where $k'$, $k''$ and $k'''$ denote the f\/irst-, second- and third-order derivatives of a real-valued func\-tion~$k$ with respect to the independent variable. The same notation will be used for vector-valued functions. Integrating twice (\ref{1:3.1}) we f\/ind
\begin{gather}\label{1:3.2}
k'^2=-\frac{8}{3}k^3+\frac{3}{2}g_2k-\frac{9}{4}g_3,
\end{gather}
where $g_2$ and $g_3$ are real constants. Every $k$ satisfying~(\ref{1:3.2}) generates the traveling wave solution
\[
u(s,t)=k\left(s-\frac{3}{4}g_2t\right)
\]
of the f\/irst equation of the hierarchy. Clearly, (\ref{1:3.2}) can be integrated in terms of the Weierstrass~$\wp$ functions, namely: if $\Delta(g_2,g_3)=-g_2^3+27g_3^2>0$, then
\[k(s)=-\frac{3}{2}\wp(s+c), \qquad s\in (2n\omega_1-c,(2n+1)\omega_1-c), \qquad n\in \Z,\]
where $\omega_1$ is the real half period and $c$ is a real constant. If $\Delta(g_2,g_3)<0$, then there are two types of solutions:
 \begin{gather*}
 k(s) = -\frac{3}{2}\wp(s+c),\qquad s\in (2n\omega_1-c,(2n+1)\omega_1-c), \qquad n\in \Z,\\
 k(s)=-\frac{3}{2}\wp(s+\omega_3+c),\qquad s\in \R,
\end{gather*}
where $\omega_1$ and $\omega_3$ are the real and the purely imaginary half periods. When $\Delta(g_2,g_3)<0$, the Weierstrass functions can be written in terms of Jacobi elliptic functions and we get
\begin{gather*}%\label{1:3.3}
k(s)=\frac{1}{2}(1+m)-\frac{3}{2}\mathrm{ns}(s+c|m)^2,\qquad
 k(s)=\frac{1}{2}\left(1-2m+3m\mathrm{cn}(s+c|m)^2\right).
\end{gather*}
 The parameter $m\in (0,1)$ is $e_3-e_2$, where $e_1>e_2>e_3$ are the three real roots of the cubic polynomial $4t^3-g_2t-g_3$. The functions of the second type are periodic, with minimal pe\-riod~$2K(m)$, where $K$ is complete elliptic integral of the f\/irst kind. The velocity of the traveling waves originated by these functions is $v_m=-(1+m(m-1))$. The f\/irst family of~\cite{Za} consists of the
traveling waves of the second type. When $m\to 0$ we obtain
\[
k(t)=\frac{1}{2}+\frac{3}{2}\mathrm{csc}(s+c)^2,
\]
and, when $m\to 1$, we f\/ind
\[k(t)=1-\frac{3}{2}\mathrm{coth}(s+c)^2,\qquad k(t)=-\frac{1}{2}+\frac{3}{2}\mathrm{sech}(s+c)^2,\]
which coincide with the solutions $(67)$ and $(69)$ of~\cite{Wa}.
Next we consider the third-order equation
\begin{gather}\label{1:3.4}
k'''+k k'=0.
\end{gather}
Again, integrating twice, we f\/ind
\begin{gather}\label{1:3.5}
k'^2=-\frac{1}{3}k^3-12g_2k-144g_3.
\end{gather}
Each function satisfying (\ref{1:3.5}) generates the traveling wave solution
\[
u(s,t)=k(s-132g_2t)
\]
of the f\/ifth-order Kaup--Kupershmidt equation. As in the previous case, the solutions of~(\ref{1:3.5}) can be expressed in terms of Weiertsrass $\wp$-functions and Jacobi elliptic functions: if $\Delta(g_2,g_3)=-g_2^3+27g_3^2>0$, then
\[k(s)=-12\wp(s+c), \qquad s\in (2n\omega_1-c,(2n+1)\omega_1-c), \qquad n\in \Z,\]
and, if $\Delta(g_2,g_3)<0$, we obtain
\begin{gather*}%\label{1:3.6}
k(s)=4(1+m)-12\mathrm{ns}(s+c|m)^2,\qquad
 k(s)=4(1-2m)+12m\mathrm{cn}(s+c|m)^2.
\end{gather*}
The velocity of the traveling waves originated by these functions is $v_m=-176(1+m(m-1))$. The second family of~\cite{Za} consists of the
traveling waves of the second type. When $m\to 0$ we obtain
\[k(t)=4-12\mathrm{csc}(s+c)^2,\]
and, when $m\to 1$, we f\/ind
\[k(t)=8-12\mathrm{coth}(s+c)^2,\qquad k(t)=-4+12\mathrm{sech}(s+c)^2,
\]
which coincide with the solutions $(68)$ and $(70)$ of~\cite{Wa}.

\section{Motion of curves in projective plane}\label{s:2}

\subsection{Curves in projective plane and their adapted frames}\label{s:2.1}

\noindent Consider a smooth parameterized curve $\gamma: I\to \RP^2$, def\/ined on some open interval $I\subset \R$ and let $G:I\to \R^3\setminus\{(0,0,0)\}$ be any lift of $\gamma$. We say that $\gamma(t)$ is an {\it inflection point} if $\mathrm{Span}(G(t),G'(t),G''(t))$ has dimension $\le 2$. From now on we will consider only curves without points of inf\/lection. Then,
\begin{gather*}%\label{2:1.1}
\Gamma = \mathrm{Det}(G,G',G'')^{-1/3}G
\end{gather*}
is the unique lift such that
\begin{gather}\label{2:1.2}
\mathrm{Det}(\Gamma,\Gamma',\Gamma'')=1.
\end{gather}
Dif\/ferentiating (\ref{2:1.2}) we see that there exist smooth functions $a,b:I\to \R$ such that
\begin{gather*}%\label{2:1.3}
\Gamma''' = a\Gamma +b\Gamma'.
\end{gather*}
The {\it projective speed} $v$ and the {\it projective arc-element} $\sigma$ are def\/ined by
\begin{gather*}%\label{2:1.4}
v=\left(a-b'/2\right)^{1/3},\qquad \sigma =vdt.
\end{gather*}
The primitives $s : I\to \R$ of the projective arc-element are the {\it projective parameters} and the zeroes of $\sigma$ are the {\it sextatic points}. A~curve without inf\/lection or sextatic points is said to be {\it generic}. Obviously, every generic curves can be parameterized by the projective parameter.

\begin{remark}\label{remark} For every point $t\in I$, there is a unique non-degenerate conic $\mathcal{C}_{t}$, the {\it osculating conic}, having fourth-order analytic contact with $\gamma$ at $\gamma(t)$. The osculating conic is def\/ined by the equation $x_1^2-2x_0x_2=0$ with respect to the homogenous coordinates of the projective frame
\[ \left(\Gamma |_{t},\Gamma '|_{t},\Gamma ''|_{t}-\frac{b(t)}{2}\Gamma|_{t}\right)\]
and, identifying $\RP^5$ with the space of plane conics, we def\/ine the osculating curve by
\[\mathcal{C}: \ t \in I \to \mathcal{C}|_{t}\in \RP^5.\]
The sextatic points are critical points of the osculating curve. The assumption on the non-existence of inf\/lection and sextatic points is rather strong from a global viewpoint. For instance, any simple closed curve in $\RP^2$ possesses f\/lex or sextatic points and a simple convex curve has at least six sextatic points \cite{TU,Um}.
\end{remark}

 The {\it canonical projective frame field} \cite{Ca} along a generic curve $\gamma$ is the $SL(3,\R)$-valued map def\/ined by
\begin{gather*}%\label{2:1.5}
F_0=v\Gamma,\qquad F_1=\frac{v'}{v}\Gamma+\Gamma',\qquad F_2=\frac{1}{2v}\left(\frac{v'^2}{v^2}-b\right)\Gamma+\frac{v'}{v^2}\Gamma+\frac{1}{v}\Gamma''.
\end{gather*}
The canonical frame is invariant with respect to changes of the parameter and
projective transformations. Furthermore, it satisf\/ies the {\it projective Frenet system}
\begin{gather}\label{2:1.6}
F'=  F\cdot
                     \begin{pmatrix}
                       0 & -k & 1 \\
                       1 & 0 & -k \\
                       0 & 1 & 0 \\
                     \end{pmatrix} v.\end{gather}
The function  $k : I\to \R$ is the {\it projective curvature}, whose explicit expression \cite{OT} is
\begin{gather*}%\label{2:1.7}
k=-\frac{1}{v} \left(S(s)+\frac{b}{2}\right),
\end{gather*}
where $s$ is a projective parameter function and $S$ is the Schwarzian derivative
\[S(f)=\left(\frac{f_{tt}}{f_t}\right)_t - \frac{1}{2}\left(\frac{f_{tt}}{f_t}\right)^2.\]

\begin{remark} The construction of the canonical frame involves only algebraic manipulations and dif\/ferentiations. So, it can be implemented in any software of symbolic calculus. In addition, the canonical frame can be constructed using the ``invariantization'' method of Fels--Olver~\cite{FO}. In other words, there is a $SL(3,\R)$-equivariant map
$\mathfrak{F}:J_h^{*}(\R ,\RP^2)\to SL(3,\R)$, def\/ined on the f\/ifth-order jet space of generic curves such that
$\mathfrak{F}\circ j^{(5)}(\gamma)$ is the projective frame along $\gamma$, for every non-degenerate~$\gamma$.
\end{remark}

\begin{figure}[t]
\centering
\includegraphics[width=4cm]{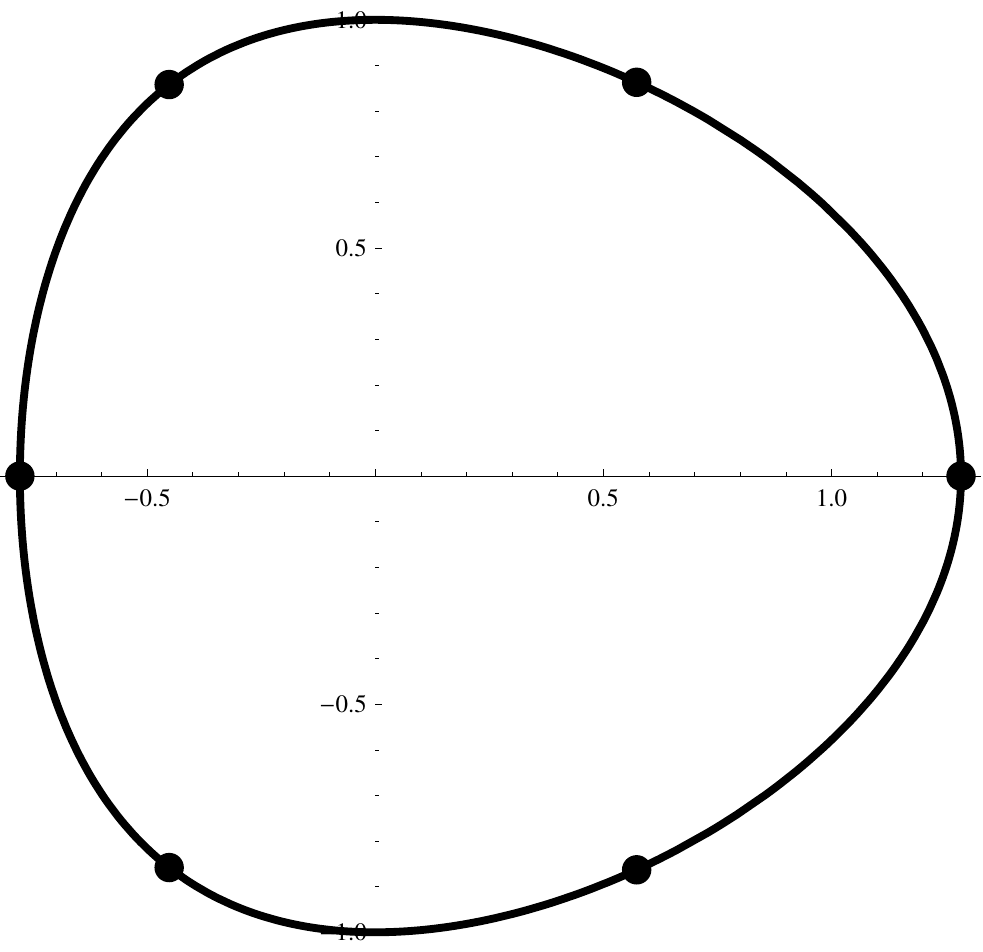}\qquad\qquad
\includegraphics[width=4cm]{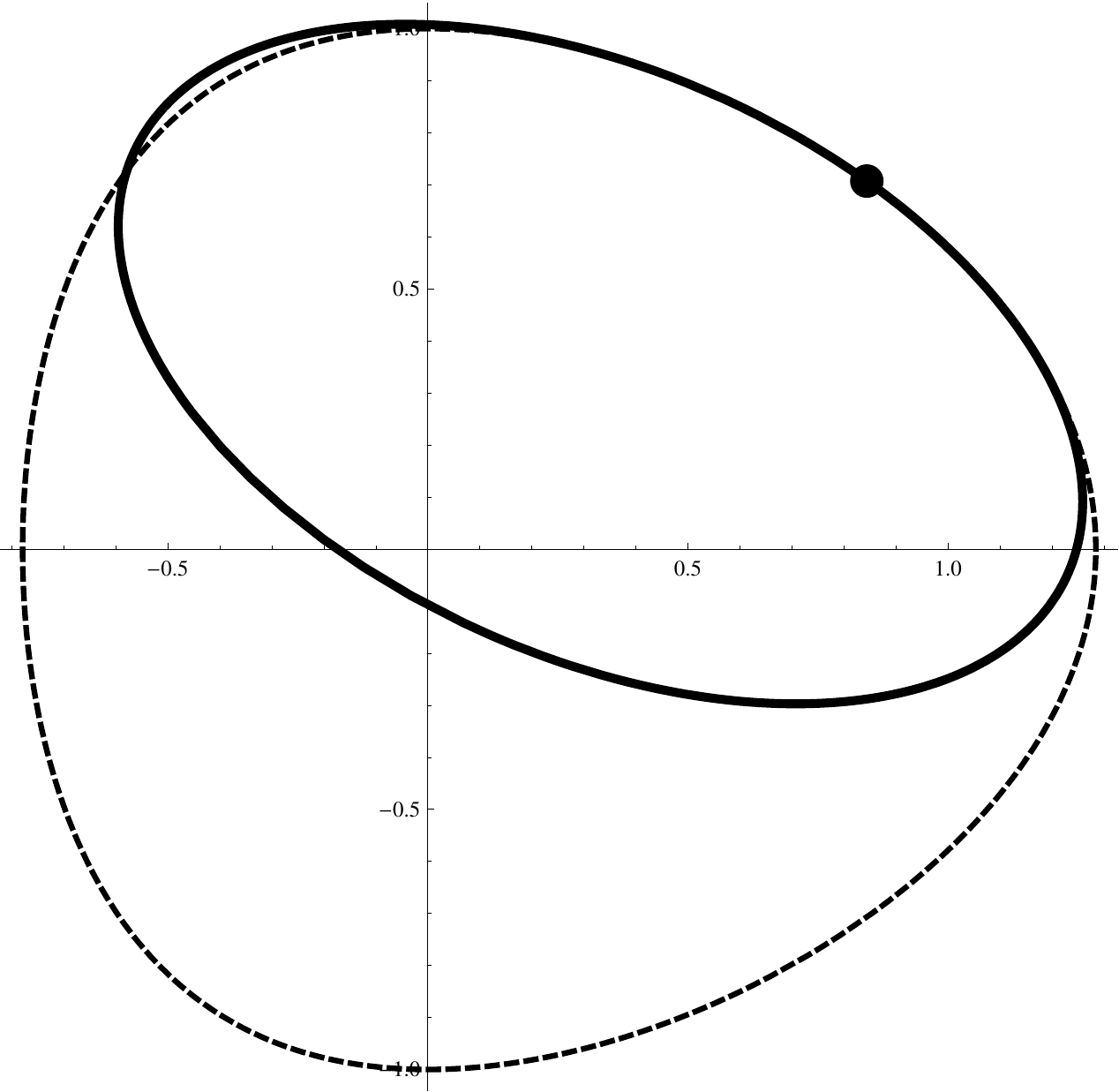}
\caption{The sextatic points and the osculating conic $\mathcal{C}(\pi/4)$.}\label{FIG1}
\end{figure}

\begin{figure}[t]
\centering
\includegraphics[width=7cm]{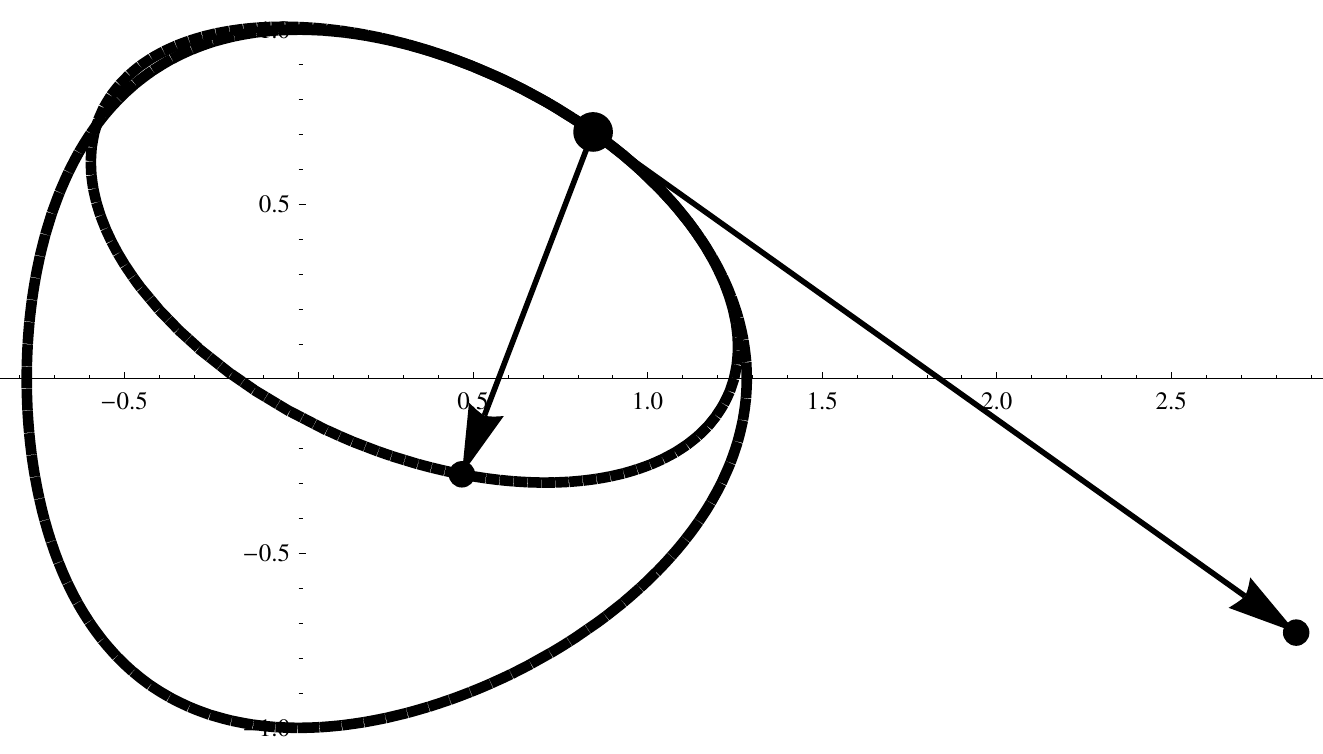}
\caption{The projective frame and the osculating conic.}\label{FIG2}
\end{figure}

\begin{figure}[t!]
\centering
\includegraphics[width=6cm]{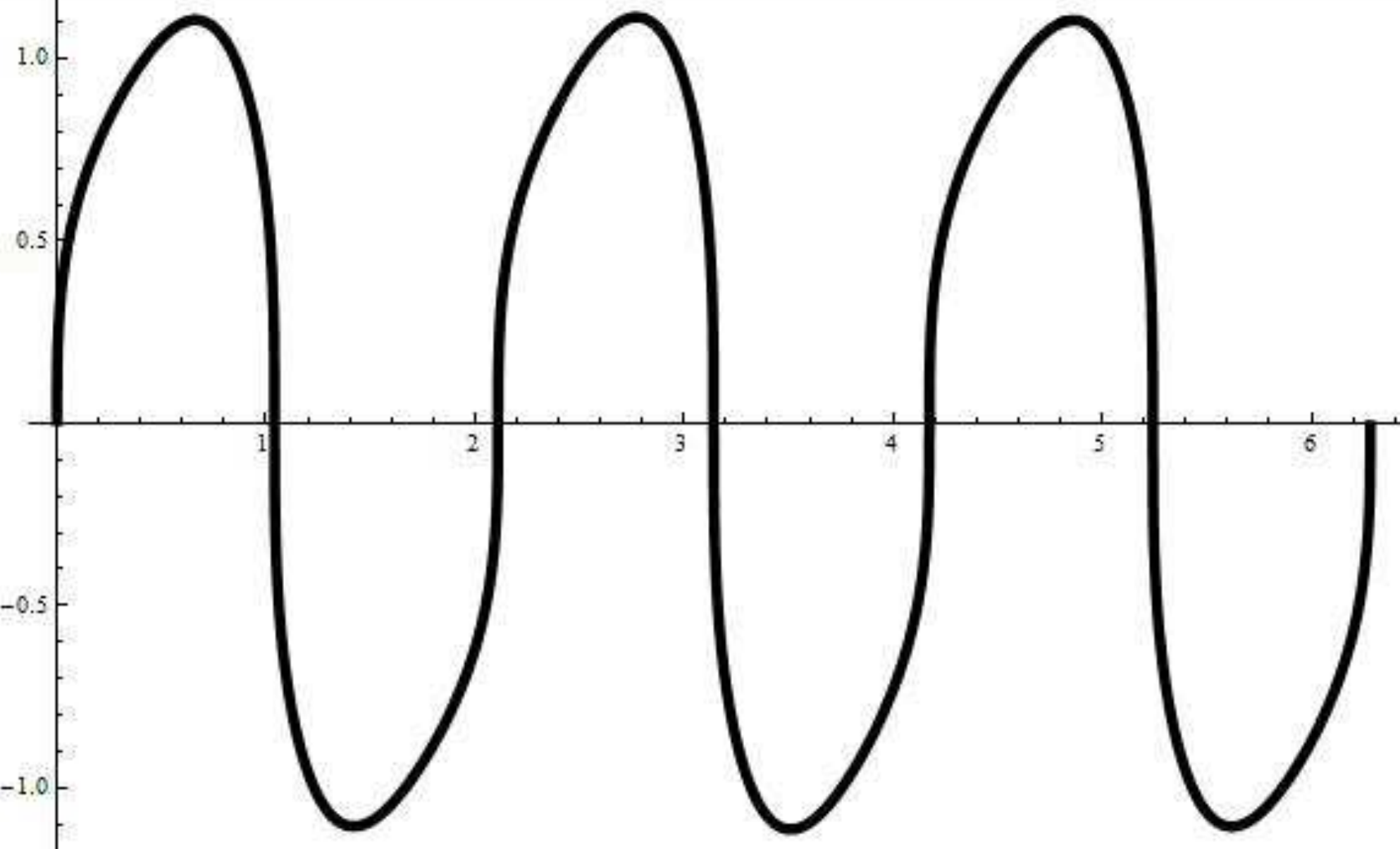}\qquad\qquad
\includegraphics[width=6cm]{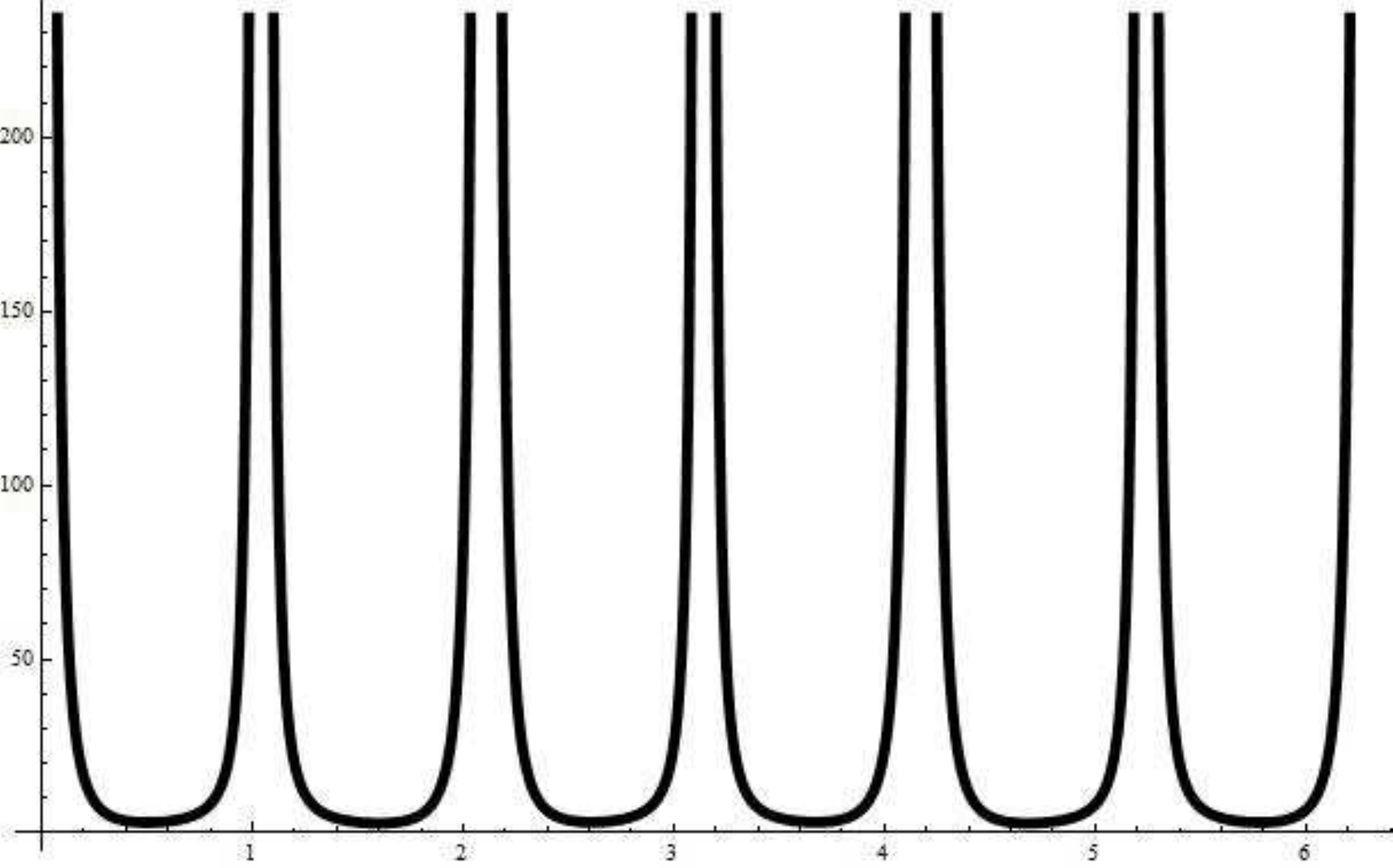}
\caption{The speed and the projective curvature.}\label{FIG3}
\end{figure}

\begin{example}
The convex simple curve
\[
\gamma : \ t\in \R \to \big[\big(\cos(t),\sin(t),e^{\frac{-\cos(t)}{4}}\big)\big]\in \RP^2
\]
has exactly six sextatic points, attained at
\begin{alignat*}{4}
 & \tau_1=0,\qquad && \tau_2\approx 1.0412803807424216,\qquad && \tau_3 \approx 2.109976014903134, &\\
 & \tau_4=\pi,\qquad && \tau_5\approx 4.173209292276453,\qquad && \tau_6\approx 5.241904926437165. &
\end{alignat*}
Figs.~\ref{FIG1} and \ref{FIG2} reproduce the curve, the sextatic points, the osculating conic and the projective frame at $t=\pi/4$. Fig.~\ref{FIG3} reproduces the projective speed and the projective curvature. The speed vanishes at the sextatic points and the curvature becomes inf\/inite at these points.
\end{example}

\begin{figure}[t]
\centering
\begin{minipage}[b]{65mm}\centering
\includegraphics[width=65mm]{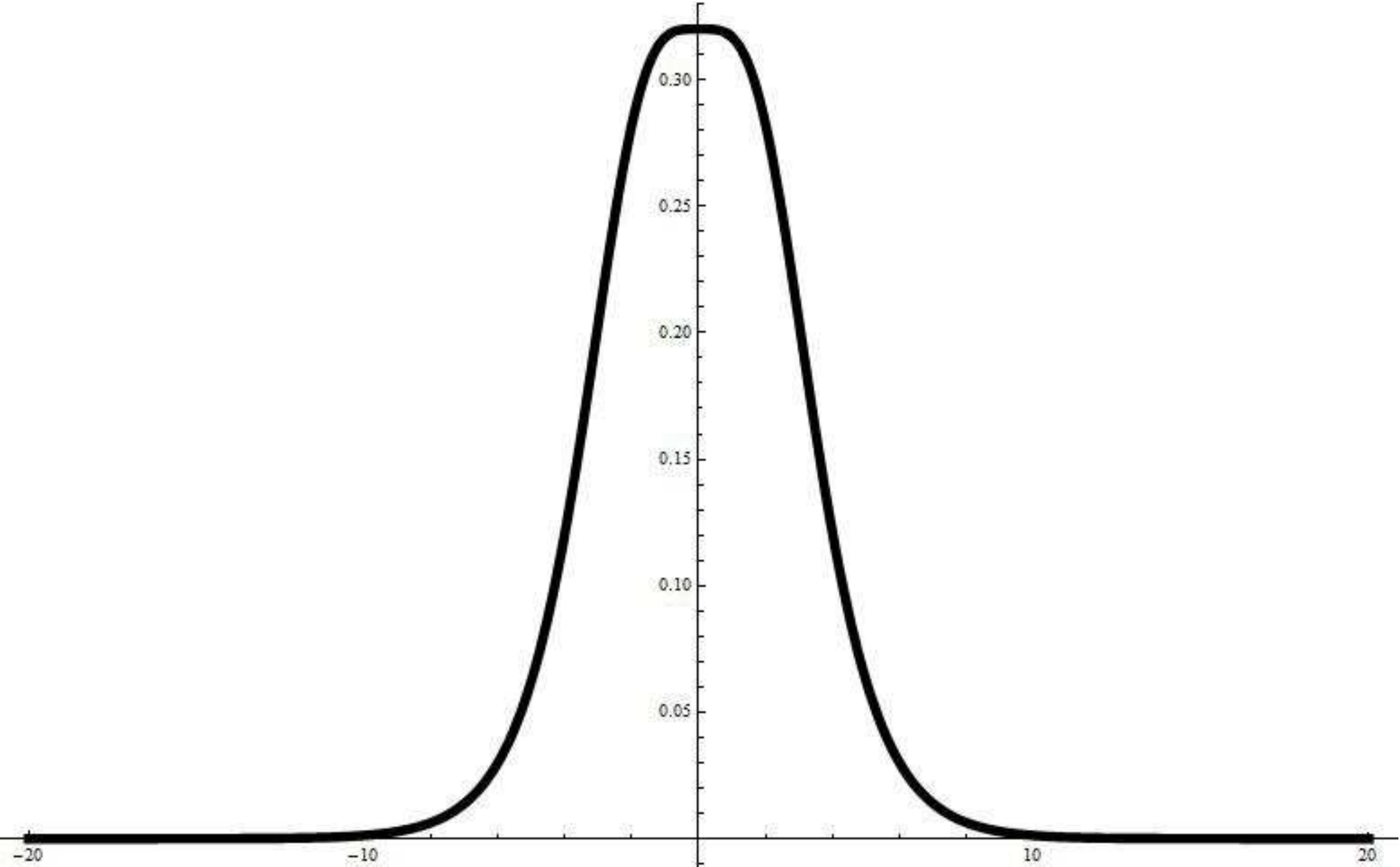}
\caption{The projective curvature.}\label{FIG4}
\end{minipage} \qquad
\begin{minipage}[b]{85mm}\centering
\includegraphics[width=4cm]{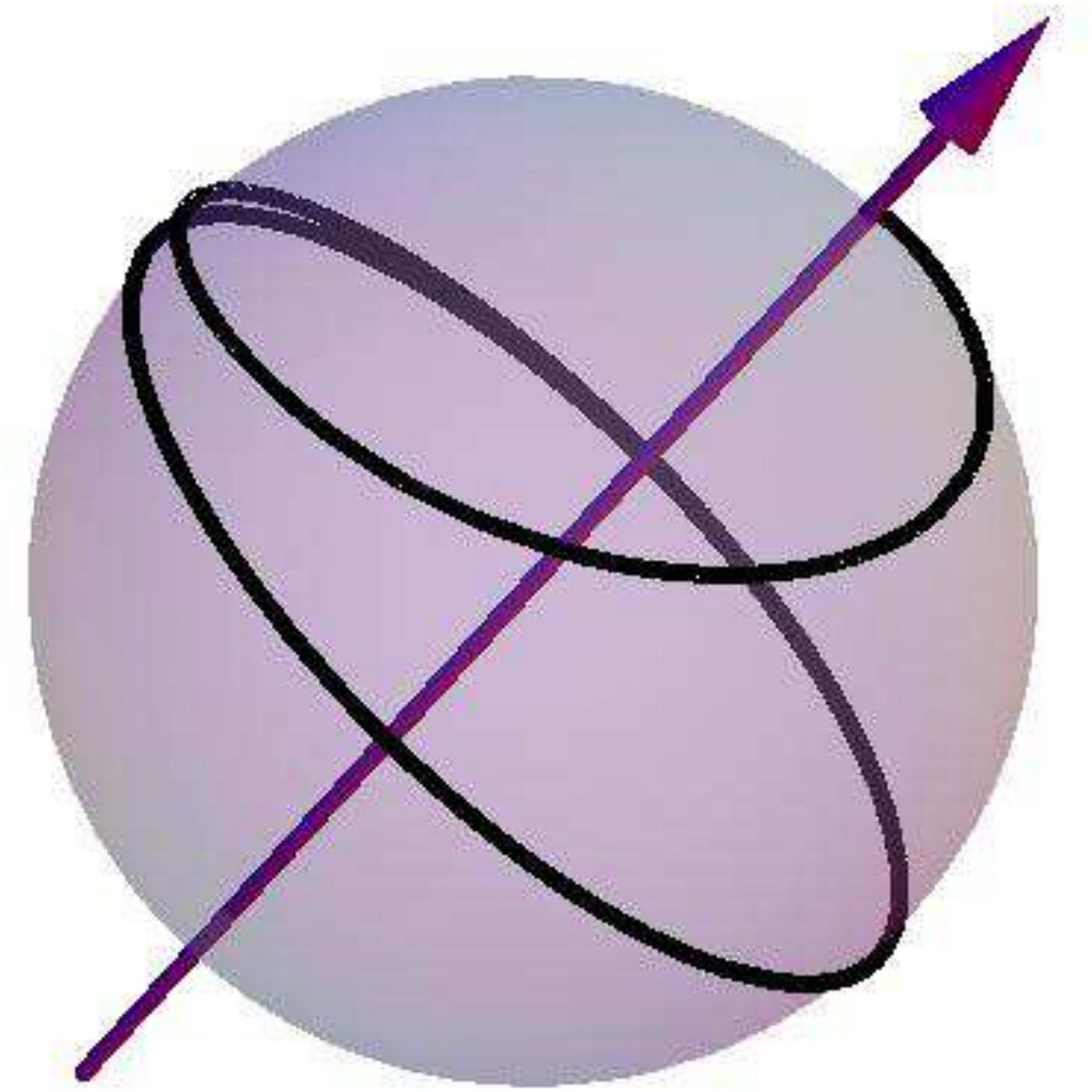}
\caption{The spherical lift of the corresponding curve.}\label{FIG5}
\end{minipage}
\end{figure}

\begin{remark} The curve is uniquely determined, up to projective congruences, by the speed and the curvature. If we assign smooth functions $v>0$ and $k$, the Frenet system~(\ref{2:1.6}) can be integrated with standard numerical routines (see Appendix~\ref{a2}). For instance, taking $v=1$ and choosing the ``anomalous'' $1$-soliton solution of the $\mathcal{K}_1$-equation~\cite{Ka,Pa}
\[k(s)=\frac{2m^2\left(1+2\cosh(m(s-m^4t)\right)}{2\left(2+\cosh(m(s-m^4t)\right)^2)},\qquad m=0.8,\]
as projective curvature (see Fig.~\ref{FIG4}), the numerical solution of the linear system~(\ref{2:1.6}) gives rise to the curve whose spherical lift is reproduced in Fig.~\ref{FIG5}.
\end{remark}

\subsection{The equation of a motion of curves}\label{s:2.2}

 A {\it motion} is
a smooth one-parameter family $\gamma(s,t)$ of projective curves such that
\[\gamma_{[t]} :\ s\in \R\to \gamma(s,t)\in \mathbb{RP}^2,\]
is generic and parameterized by the projective parameter, for every $t\in I$.
Denoting by $F_{[t]}:\R\to SL(3,\R)$ and $k_{[t]}:\R\to \R$ the projective frame and the projective curvature of $\gamma_{[t]}$ we consider the {\it projective frame} and the {\it projective curvature} of the motion, def\/ined by
\begin{gather*}%\label{2:2.1}
 \mathcal{F}: \ (s,t)\in \R\times I\to F_{[t]}(s)\in SL(3,\R)
 \end{gather*}
and
\begin{gather*}%\label{2:2.2}
\kappa: \ (s,t)\in \R\times I\to k_{[t]}(s)\in \R.
\end{gather*}
The projective frame satisf\/ies
\begin{gather}\label{2:2.3}
\mathcal{F}^{-1}d\mathcal{F} = \mathcal{K}(s,t)ds + \Phi(s,t) dt,
\end{gather}
where
\begin{gather}\label{2:2.4}
\mathcal{K}=\begin{pmatrix}
        0 & -\kappa & 1 \\
        1 & 0 & -\kappa \\
        0 & 1 & 0
      \end{pmatrix},\qquad
    \Phi=\begin{pmatrix}
             \phi^0_0 & \phi^0_1 & \phi^0_2 \\
             \phi^1_0 & \phi^1_1 & \phi^1_2 \\
             \phi^2_0 & \phi^2_1 & -(\phi^0_0+\phi^1_1)
           \end{pmatrix}.
\end{gather}
The coef\/f\/icient $\phi^2_0$ is said to be the {\it normal velocity of the motion} and it will be denoted by~$\upsilon$.

\begin{proposition}\label{prop1}
The curvature of a motion of projective curves with normal velocity $\upsilon$ satisfies
\begin{gather}\label{2:2.5}
\partial_t\kappa = \Theta(\upsilon,\kappa)+\lambda \kappa_s,
\end{gather}
where $\lambda$ is a real constant, the {\it internal parameter}. Conversely, if $\kappa$ is a solution of~\eqref{2:2.5} then there is a motion $\gamma$ with normal speed $\upsilon$, internal parameter $\lambda$ and curvature~$\kappa$. Moreover,~$\gamma$ is unique up to projective transformations.
\end{proposition}

\begin{proof}
Dif\/ferentiating (\ref{2:2.3}) we obtain
\begin{gather}\label{2:2.6}
\partial_s\Phi-\partial_t \mathcal{K}+[\mathcal{K},\Phi]=0,
\end{gather}
which implies
\begin{gather}
(\phi^0_0)_s -\kappa \phi^1_0 + \phi^2_0 - \phi^0_1 =  0,\nonumber\\
(\phi^0_2)_s-2\phi^0_0-\phi^1_1+\kappa(\phi^0_1-\phi^1_2)=0,\nonumber\\
(\phi^1_0)_s+\phi^0_0-\kappa \phi^2_0-\phi^1_1=0,\nonumber\\
(\phi^1_1)_s+\phi^0_1+\kappa(\phi^1_0-\phi^2_1)-\phi^1_2=0,\nonumber\\
(\phi^2_0)_s-\phi^2_1=0,\nonumber\\
(\phi^2_1)_s+\phi^0_0+\kappa \phi^2_0+2\phi^1_1=0,\nonumber\\
(\phi^0_1)_s-(\phi^1_2)_s-3\kappa\phi^1_1+\phi^2_1+\phi^1_0-2\phi^0_2=0,\label{2:2.7}
 \end{gather}
and
\begin{gather}\label{2:2.8}
\partial_t \kappa+(\phi^1_2)_s-\phi^1_0+\kappa(\phi^0_0+2\phi^1_1)+\phi^0_2=0.
\end{gather}
If we set $\upsilon=\phi^2_0$, then (\ref{2:2.7}) gives
\begin{gather}
\phi^0_0 =\left(\frac{1}{3}\kappa+\frac{8}{9}\kappa\kappa_s+\frac{1}{9}\kappa_{3s}\right)\upsilon+
\left(\frac{1}{2}\kappa_s+\frac{8}{9}\kappa^2\right)\upsilon_s+\left(\frac{1}{6}
+\frac{5}{6}\kappa_s\right)\upsilon_{2s}+\frac{5}{9}\kappa\upsilon_{3s}+\frac{1}{18}\upsilon_{5s},
\nonumber\\
\phi^1_0 =\lambda -\frac{1}{9}\int_0^s \big(\kappa\upsilon_{3s}+4\kappa^2\upsilon_s\big)ds-\left(\frac{1}{9}\kappa_{2s}+\frac{4}{9}\kappa^2\right)\upsilon-
\left(\frac{1}{2}+\frac{7}{18}\kappa_s\right)\upsilon_s
-\frac{4}{9}\kappa\upsilon_{2s}-\frac{1}{18}\upsilon_{4s},\nonumber\\
\phi^0_1  = -\lambda\kappa+\frac{\kappa}{9}\int_0^s \big(\kappa\upsilon_{3s}+4\kappa^2\upsilon_s\big)ds
+\left(1+\frac{4}{9}\kappa^3+\frac{1}{3}\kappa_s+\frac{8}{9}\kappa_s^2+\kappa\kappa_{2s}+\frac{1}{9}\kappa_{4s}\right)\upsilon \label{2:2.9}\\
\hphantom{\phi^0_1  =}{} + \left(\frac{5}{6}\kappa+\frac{55}{18}\kappa\kappa_s+\frac{11}{8}\kappa_{3s}\right)\upsilon_s\!+\frac{4}{3}\big(\kappa^2+\kappa_{2s}\big)\upsilon_{2s}\! +
\left(\frac{1}{6}+\frac{25}{18}\kappa_s\right)\upsilon_{3s}\!+\frac{11}{8}\kappa\upsilon_{4s}\!+\frac{1}{18}\upsilon_{6s},\nonumber
\end{gather}
and
\begin{gather}
\phi^1_1 =-\frac{2}{3}\kappa\upsilon-\frac{1}{3}\upsilon_{2s},\nonumber\\
\phi^2_1 =\lambda-\frac{1}{9}\int_0^s \big(\kappa\upsilon_{3s}+4\kappa^2\upsilon_s\big)ds-\left(\frac{4}{9}\kappa^2+\frac{1}{9}\kappa_{2s}\right)\upsilon
+\left(\frac{1}{2}-\frac{7}{18}\kappa_s\right)\upsilon_s
-\frac{4}{9}\kappa\upsilon_{2s}
-\frac{1}{18}\upsilon_{4s},\nonumber\\
\phi^0_2 =\lambda -\frac{1}{9}\int_0^s \big(\kappa\upsilon_{3s}+4\kappa^2\upsilon_s\big)ds+\left(\frac{5}{9}\kappa^2+\frac{2}{9}\kappa_{2s}\right)\upsilon+\frac{7}{9}\kappa_s\upsilon_s+
\frac{8}{9}\kappa\upsilon_{2s}+\frac{1}{9}\upsilon_{4s},\nonumber\\
\phi^1_2 =-\lambda\kappa+\frac{\kappa}{9}\int_0^s \big(\kappa\upsilon_{3s}+4\kappa^2\upsilon_s\big)ds+\left(1+\frac{4}{9}\kappa^3-\frac{1}{3}\kappa_s+\frac{8}{9}\kappa_s^2+\kappa\kappa_{2s}
+\frac{1}{9}\kappa_{4s}\right)\upsilon \label{2:2.10}\\
\phantom{\phi^1_2 =}{}
+\left(\!{-}\frac{5}{6}\kappa+\frac{55}{18}\kappa\kappa_s+\frac{11}{18}\kappa_{3s}\right)\!\upsilon_s\!+\frac{4}{3}\big(\kappa^2\!+\kappa_{2s}\big)\upsilon_{2s}\!-
\left(\frac{1}{6}-\frac{25}{18}\kappa_s\right)\!\upsilon_{3s}\!+\frac{11}{18}\kappa\upsilon_{4s}\!+\frac{1}{18}\upsilon_{6s},\nonumber
\end{gather}
where $\lambda$ is a real constant. From (\ref{2:2.9}) and (\ref{2:2.10}) we deduce that (\ref{2:2.8}) is satisf\/ied if and only if
\[
 \partial_t\kappa   =\Theta(\upsilon,\kappa)+\lambda\kappa_s.
\]
Conversely, if $\kappa$ is a solution of~(\ref{2:2.5}) and if we def\/ine $\phi^i_j$, $\mathcal{K}$ and $\Phi$ as in~(\ref{2:2.4}), (\ref{2:2.9}) and (\ref{2:2.10}), then~$K$ and~$\Phi$ satisfy~(\ref{2:2.6}). Using Frobenius theorem we deduce the existence of a smooth map
\[\mathcal{F}: \ \R\times I\to SL(3,\R)\]
such that
$\mathcal{F}^{-1}d\mathcal{F}=\mathcal{K}ds+\Phi dt$. The map~$\mathcal{F}$ is unique up to left multiplication by an element of~$SL(3,\mathbb{R})$. Setting $\gamma(s,t)=[F_0(s,t)]$ we have a motion of projective curves with curvature~$\kappa$, normal velocity~$\upsilon$, internal parameter~$\lambda$ and projective frame $\mathcal{F}$. This yields the required result.
 \end{proof}

\subsection{Local motions}\label{s:2.3}

From the proof of Proposition \ref{prop1} we see that the $\Phi$-matrix of a motion of projective curves can be written as
\begin{gather*}%\label{2:2.11}
\Phi=\widetilde{\Phi}(\upsilon,\kappa)+\lambda\mathcal{K}(\kappa),
\end{gather*}
where the coef\/f\/icients of $\widetilde{\Phi}(\upsilon,\kappa)$ are the integro-dif\/ferential operators
\begin{gather*}
\tilde{\phi}^2_0(\upsilon,\kappa) =\upsilon,\\
\tilde{\phi}^0_0(\upsilon,\kappa) =\left(\frac{1}{3}\kappa+\frac{8}{9}\kappa\kappa_s+\frac{1}{9}\kappa_{3s}\right)\upsilon+
\left(\frac{1}{2}\kappa_s+\frac{8}{9}\kappa^2\right)\upsilon_s\!
+\left(\frac{1}{6}+\frac{5}{6}\kappa_s\right)\upsilon_{2s}\!+\frac{5}{9}\kappa\upsilon_{3s}\!+\frac{1}{18}\upsilon_{5s},\\
\tilde{\phi}^1_0(\upsilon,\kappa) = -\frac{1}{9}\!\int_0^s \! \big(\kappa\upsilon_{3s}\!+4\kappa^2\upsilon_s\big)ds-\left(\frac{1}{9}\kappa_{2s}+\frac{4}{9}\kappa^2\right)\upsilon
-\left(\frac{1}{2}+\frac{7}{18}\kappa_s\right)\upsilon_s\!
-\frac{4}{9}\kappa\upsilon_{2s}\!-\frac{1}{18}\upsilon_{4s},\\
 \tilde{\phi}_1^0(\upsilon,\kappa) = \frac{\kappa}{9}\int_0^s \big(\kappa\upsilon_{3s}+4\kappa^2\upsilon_s\big)ds+
 \left(1+\frac{4}{9}\kappa^3+\frac{1}{3}\kappa_s+\frac{8}{9}\kappa_s^2+\kappa\kappa_{2s}+\frac{1}{9}\kappa_{4s}\right)\upsilon \\
\phantom{\tilde{\phi}_1^0(\upsilon,\kappa) =}{}
+\left(\frac{5}{6}\kappa+\frac{55}{18}\kappa\kappa_s+\frac{11}{8}\kappa_{3s}\right)\upsilon_s+\frac{4}{3}\big(\kappa^2+\kappa_{2s}\big)\upsilon_{2s}+
\left(\frac{1}{6}+\frac{25}{18}\kappa_s\right)\upsilon_{3s}\\
\phantom{\tilde{\phi}_1^0(\upsilon,\kappa) =}{}
+\frac{11}{8}\kappa\upsilon_{4s}+\frac{1}{18}\upsilon_{6s},\\
\tilde{\phi}^1_1(\upsilon,\kappa) =-\frac{2}{3}\kappa\upsilon-\frac{1}{3}\upsilon_{2s},\\
\tilde{\phi}^2_1(\upsilon,\kappa) =-\frac{1}{9}\!\int_0^s\! \big(\kappa\upsilon_{3s}\!+4\kappa^2\upsilon_s\big)ds-
\left(\frac{4}{9}\kappa^2\!+\frac{1}{9}\kappa_{2s}\right)\upsilon+\left(\frac{1}{2}-\frac{7}{18}\kappa_s\right)\upsilon_s\!
-\frac{4}{9}\kappa\upsilon_{2s}\!
-\frac{1}{18}\upsilon_{4s},\\
\tilde{\phi}^0_2(\upsilon,\kappa) = -\frac{1}{9}\int_0^s \big(\kappa\upsilon_{3s}+4\kappa^2\upsilon_s\big)ds+\left(\frac{5}{9}\kappa^2+\frac{2}{9}\kappa_{2s}\right)\upsilon+\frac{7}{9}\kappa_s\upsilon_s+
\frac{8}{9}\kappa\upsilon_{2s}+\frac{1}{9}\upsilon_{4s},\\
\tilde{\phi}^1_2(\upsilon,\kappa) =\frac{\kappa}{9}\int_0^s \big(\kappa\upsilon_{3s}+4\kappa^2\upsilon_s\big)ds
+\left(1+\frac{4}{9}\kappa^3-\frac{1}{3}\kappa_s+\frac{8}{9}\kappa_s^2+\kappa\kappa_{2s}+\frac{1}{9}\kappa_{4s}\right)\upsilon \\
\phantom{\tilde{\phi}^1_2(\upsilon,\kappa) =}{}
+\left(-\frac{5}{6}\kappa+\frac{55}{18}\kappa\kappa_s+\frac{11}{18}\kappa_{3s}\right)\upsilon_s
+\frac{4}{3}\big(\kappa^2+\kappa_{2s}\big)\upsilon_{2s}-
\left(\frac{1}{6}-\frac{25}{18}\kappa_s\right)\upsilon_{3s}\\
\phantom{\tilde{\phi}^1_2(\upsilon,\kappa) =}{}
+\frac{11}{18}\kappa\upsilon_{4s}+\frac{1}{18}\upsilon_{6s}.
\end{gather*}
Bearing in mind the def\/inition (\ref{P[u]}) of the linear subspace $P[\mathbf{u}]$, we deduce the existence of linear operators
$\mathcal{M}^i_j:P[\mathbf{u}]\to J[\mathbf{u}]$ such that
\begin{gather*}%\label{2:3.1}
\widetilde{\phi}^i_j(\mathfrak{p}|_{j_s(\kappa)},\kappa)=\mathcal{M}^i_j(\mathfrak{p})|_{j_s(\kappa)},
\end{gather*}
for every $\mathfrak{p}\in P[\mathbf{u}]$. This implies the following corollary.

\begin{corollary}
If $\mathfrak{p}$ belongs to $P[\mathbf{u}]$ and if $\kappa$ is solution of the evolution equation
\begin{gather*}%\label{2:3.2}
\partial_t\kappa = \mathcal{S}(\mathfrak{p})|_{j_s(\kappa)}+\lambda \kappa_s
\end{gather*}
then, there is a motion $\gamma$, uniquely defined up to projective transformations, with curvature $\kappa$ and normal velocity $\mathfrak{p}|_{j_s(\kappa)}$. Motions of this type are said to be {\it local}.
\end{corollary}

\begin{remark} Local motions are the integral curves of {\it local vector fields} on the inf\/inite-dimensional space $\mathcal{P}$ of unit-speed generic curves of $\RP^2$. More precisely, if we take any $\mathfrak{p}\in P[\mathbf{u}]$ and any real constant $\lambda$ then there is a unique vector f\/ield $X_{\mathfrak{p},\lambda}$ on $\mathcal{P}$ whose integral curve through $\gamma_{[0]}\in \mathcal{P}$ is the local motion $\gamma$ such that:
\begin{itemize}\itemsep=0pt
\item $\gamma(s,0)=\gamma_{[0]}(s)$;
\item its curvature $\kappa$ is the solution of the Cauchy problem
\[
\kappa_t+\mathcal{S}(\mathfrak{p})|_{j_s(\kappa)}+\lambda\kappa_s = 0,\qquad
\kappa(s,0)=k_{[0]}(s);
\]
\item its normal speed is $\mathfrak{p}|_{j_s(\kappa)}$.
\end{itemize}
\end{remark}

\begin{definition} We say that  $X_{\mathfrak{p},\lambda}$  is the {\it local vector field} with potential $\mathfrak{p}$ and spectral parame\-ter~$\lambda$. The dynamics of a local vector f\/ield is governed by the {\it induced evolution equation}
\[
\kappa_t+\mathcal{S}(\mathfrak{p})|_{j_s(\kappa)}+\lambda\kappa_s = 0.
\]
\end{definition}

From these observations and using Proposition \ref{proposition1} we have the following result.

\begin{theorem}
For every $n\in \mathbb{N}$ the local vector field $X_{\mathfrak{v}_n,\lambda_n}$ defined by the polynomial differential function $\mathfrak{v}_n\in P[\mathbf{u}]$ and by the spectral parameter $\lambda_n$ induces the $n$-th equation of the Kaup--Kupershmidt hierarchy.
\end{theorem}

\section{Congruence motions}
Consider a local dynamics with potential $\mathfrak{p}$ and internal parameter $\lambda$. A curve $\tilde{\gamma}$ which evolves without changing its shape (by projective transformations) is said to be a {\it congruence curve} of the f\/low. Denote by $\tilde{k}$ the curvature of $\tilde{\gamma}$ and by $\kappa(s,t)$ the curvature of the evolution~$\gamma(s,t)$ of~$\tilde{\gamma}(s)$. If $\tilde{k}$ is non constant, then
\[
\kappa(s,t)=\tilde{k}(s + vt),
\]
for some constant $v$. So, $\kappa$ is a traveling wave solution of the induced evolution equation and $\tilde{k}$ satisf\/ies the ordinary dif\/ferential equation
\begin{gather}\label{3.1}
\Theta(\mathfrak{p}|_{j_s(u)},u)+(\lambda+v) u_s = 0.
\end{gather}
Unit-speed generic curves whose curvature satisf\/ies (\ref{1:3.1}) or (\ref{1:3.4}) are examples of congruence curves of the f\/irst f\/low of the hierarchy. On the other hand, (\ref{1:3.1}) is the Euler--Lagrange equation of the invariant functional def\/ined by the integral of the projective arc-element $\sigma$ \cite{Ca,MN1}. This implies the following corollary.

\begin{corollary}
Every critical curve of the functional
\[\gamma \to \int_{\gamma}\sigma\]
is a congruence curve of the first flow of the Kaup--Kupershmidt hierarchy.
\end{corollary}

The projective frame $F(s,t)$ satisf\/ies
\begin{gather}\label{3.2}
F^{-1}dF = \tilde{K}(s+vt)dt+\tilde{\Phi}(s+vt)ds,
\end{gather}
where the $\mathfrak{sl}(3,\R)$-valued functions $\tilde{K}$ and $\tilde{\Phi}$ are def\/ined as in~(\ref{2:2.4}), (\ref{2:2.9}) and (\ref{2:2.10}), with normal speed $\mathfrak{p}|_{j_s(\tilde{k})}$.
We def\/ine the {\it Hamiltonian} by
\begin{gather*}%\label{3.3}
 H=\tilde{\Phi}-v\tilde{K}:\ \R\to \mathfrak{sl}(3,\R).
 \end{gather*}
The integrability condition of (\ref{3.2}) is the Lax equation
\begin{gather*}%\label{3.4}
H'=[H,\tilde{K}]
\end{gather*}
which implies the {\it conservation law}
\begin{gather*}%\label{3.5}
\tilde{F}\cdot H \cdot \tilde{F}^{-1}=\xi,
\end{gather*}
where $\tilde{F}$ is the projective frame of $\tilde{\gamma}$ and
 $\xi$ is a f\/ixed element of $\mathfrak{sl}(3,\R)$, the {\it momentum}
 of the congruence curve $\tilde{\gamma}$. In particular,  $H$ and $\xi$ have the same spectrum. From now on we assume that $F(0)={\rm Id}_{3\times 3}$.

\begin{proposition} The motion of a congruence curve $\tilde{\gamma}$ is given by
\begin{gather}\label{3.6}
\gamma(s,t)=\mathrm{Exp}(t\xi)\cdot \tilde{\gamma}(s+vt).
\end{gather}
\end{proposition}

\begin{proof} Def\/ine $\gamma(s,t)$ as in (\ref{3.6}) and set
\[
F(s,t) = \mathrm{Exp}(t\xi)\cdot \widetilde{F}(s+vt) \qquad \forall\, (s,t)\in \R\times I.
\]
Since $F$ is a lift of $\gamma(s,t)$, it suf\/f\/ices to prove that $F$ satisf\/ies (\ref{3.2}). From the def\/inition we deduce
\[
F^{-1}\partial_sF|_{(s,t)}=\tilde{K}(s+vt)
\]
and
\begin{gather*}
F\partial_tF|_{(s,t)}  =(\tilde{F}(s+vt)^{-1}\cdot \mathrm{Exp}(-t\xi))\cdot(\mathrm{Exp}(t\xi)\xi\cdot
\widetilde{F}(s+vt)+v\mathrm{Exp}(t\xi)\cdot\partial_s\widetilde{F}|_{s+vt}) \\
\phantom{F\partial_tF|_{(s,t)}}{}   = \tilde{F}(s+vt)^{-1}\cdot \xi\cdot \widetilde{F}(s+vt)+v\tilde{K}(s+vt) = H(s+vt)+v\tilde{K}(s+vt)\\
\phantom{F\partial_tF|_{(s,t)}}{}
 = \tilde{\Phi}(s+vt).
\end{gather*}
This implies the required result.
\end{proof}

 We now prove the following proposition.

\begin{proposition} If $\tilde{k}$ is a non-constant real-analytic solution of \eqref{3.1} and if $\xi$ has three distinct eigenvalues then the corresponding congruence curve can be found by quadratures.
\end{proposition}
\begin{proof}
It suf\/f\/ices to show that the solution of the linear system
\[\tilde{F}_s=\tilde{F}\cdot \tilde{K},\qquad \tilde{F}(0)={\rm Id}_{3\times 3}\]
can be constructed from $\tilde{k}$ by algebraic manipulations, dif\/ferentiations and integrations of functions involving $\tilde{k}$ and its derivatives. This can be shown with the following reasoning:

  {\it Step I.} The Hamiltonian $H$ can be directly constructed from the prolongation $j(\tilde{k})$ of the curvature and from the potential~$\mathfrak{p}$. Then, we compute the momentum $\xi = H(0)$ and its eigenvalues~$\tau_0$,~$\tau_1$ and $\tau_2$. For each $\tau_j$ we choose row vectors $H^{j_1}$ and $H^{j_2}$ of $H$ such that{\samepage
\[
S_j=\big({}^tH^{j_1}-\tau_j \epsilon_{j_1}\big)\times \big(^tH^{j_2}-\tau_j\epsilon_{j_2}\big)\neq 0,
\]
where $(\epsilon_0,\epsilon_1,\epsilon_2)$ is the standard basis of $\mathbb{C}^3$.}

 {\it Step II.} Next, we def\/ine the $S$-matrix
\[
S=(S_1,S_2,S_3) : \  I \to \mathfrak{gl}(3,\C).
\]
This is a real-analytic map which can be computed in terms of $\tilde{k}$ and its derivatives. Subsequently, we def\/ine the $\Sigma$-matrix
\begin{gather}\label{3.7}
\Sigma = F\cdot S : \ \R\to \mathfrak{gl}(3,\C).
\end{gather}
The columns $\Sigma_j(s)$ are eigenvectors of the momentum $\xi$, with eigenvalue $\tau_j$, for each $j=0,1,2$ and every $s$. In particular,  $\Sigma_j$ has constant direction. Hence, there exist real-analytic complex valued functions $r_j$ such that
\begin{gather}\label{3.8}
\Sigma_j'=r_j\Sigma_j,\qquad j=0,1,2.
\end{gather}
Dif\/ferentiating $\Sigma = F\cdot S$ and using the structure equations satisf\/ied by $F$ we deduce
\begin{gather*}%\label{3.9}
S' + \tilde{K}\cdot S = S\cdot \Delta(r_0,r_1,r_2),
\end{gather*}
where $\Delta(r_0,r_1,r_2)$ is the diagonal matrix with elements $r_0$, $r_1$ and $r_2$. This shows that $r_0$, $r_1$ and $r_2$ can be computed in terms of $\tilde{k}$ and its derivatives.

  {\it Step III.} We compute the {\it integrating factors}
\begin{gather*}%\label{3.10}
\rho_j(s) =\mathrm{Exp}\left( \int_0^s r_j(u)du\right),\qquad j=0,1,2.
\end{gather*}
Equation (\ref{3.8}) implies
\[\Sigma_j=\rho_jC_j,\qquad j=0,1,2,\]
where $C_0$, $C_1$ and $C_3$ are constant vectors of $\C^3$. We then have
\begin{gather}\label{3.11}
\Sigma=C\cdot \Delta(\rho_0,\rho_1,\rho_2),
\end{gather}
where $C$ is a f\/ixed element of $GL(3,\C)$. Substituting (\ref{3.7}) into (\ref{3.11}) we obtain
\begin{gather*}%\label{3.12}
F = M(0)^{-1}\cdot M,
\end{gather*}
where the $M$-matrix is def\/ined by
\begin{gather*}%\label{3.13}
M\cdot S= \Delta(\rho_0,\rho_1,\rho_2).
\end{gather*}
All the steps involve only linear algebra manipulations, dif\/ferentiations and the quadratures of the functions $r_0$, $r_1$ and $r_2$, as claimed.
 \end{proof}

\begin{example} We illustrate the integration of the congruence curves with projective curvature
\[
k_m(s)=\frac{1}{2}\left(1-2m+3m\mathrm{cn}(s+c|m)^2\right),\qquad  m\in (0,1).
\]
Computing the $H$-matrix we obtain
\[H=
      \begin{pmatrix}
        h_1^1 & h^1_2 & h^1_3 \\
        0 & -2h^1_1 & h^2_1 \\
        9 & 0 & h_1^1 \\
      \end{pmatrix},\]
where the coef\/f\/icients are given by
\begin{gather*}
h_1^1  =  \frac{3}{2} \left(3 \text{dn}(s|m)^2+m-2\right),\qquad
h^1_2 =9-9 m \text{cn}(s|m) \text{dn}(s|m) \text{sn}(s|m),\\
h^3_1 =\frac{9}{4} \left(-3
\text{dn}(s|m)^4-2 (m-2) \text{dn}( s|m)^2+m^2\right),\\
h^2_1 =9 m \text{cn}( s|m) \text{dn}(s|m) \text{sn}(s|m)+9.
 \end{gather*}
Next we compute the momentum and we get
\[\xi=
\begin{pmatrix}
 \frac{3}{2} (m+1) & 9 & \frac{9}{4} (m-1)^2 \\
 0 & -3 (m+1) & 9 \\
 9 & 0 & \frac{3}{2} (m+1)
\end{pmatrix}.
\]
The discriminant of its characteristic polynomial is
\[
\delta=-31 + m (6 + m (7 + (-6 + m) m)).
\]
 From now on we assume $\delta\neq 0$. If $\delta<0$ the momentum has one real eigenvalue and two complex conjugate eigenvalues, otherwise the momentum has three distinct real eigenvalues. The eigenvalues are:
 \begin{gather*}
\tau_0 =-3\frac{\sqrt[3]{2}(n_1+n_2)^{2/3}+2n_3}{\sqrt[3]{4}(n_1+n_2)^{1/3}},\qquad
\tau_1 =3\frac{\sqrt[3]{2}(1-i\sqrt{3})(n_1+n_2)^{2/3}+2(1+i\sqrt{3})n_3}{2\sqrt[3]{4}(n_1+n_2)^{1/3}},\\
\tau_2 =3\frac{\sqrt[3]{2}(1+i\sqrt{3})(n_1+n_2)^{2/3}+2(1-i\sqrt{3})n_3}{2\sqrt[3]{4}(n_1+n_2)^{1/3}},
\end{gather*}
where
 \begin{gather*}
n_1  =3\left(\sqrt{-3 (-31 + m (6 + m (7 + (-6 + m) m)))}-9\right),\\
n_2  =(2-m) (m+1) (2 m-1),\qquad
n_3  =((m-1) m+1).
\end{gather*}
We set
\[
S_j=(^tH^2-\tau_j\epsilon_2)\times (^tH^3-\tau_j\epsilon_3),\qquad j=0,1,2,
\]
so that
\begin{gather*}
S_j^1 =-\frac{1}{2} \left(9 \text{dn}( s|m)^2+3 (m-2)-2 \tau_j \right) \left(9 \text{dn}( s|m)^2+3 (m-2)+\tau_j \right),\\
S_j^2 =81
\left(m  \text{cn}(s|m) \text{dn}(s|m) \text{sn}(s|m)+1\right),\qquad
S_j^3 =9 \left(9  \text{dn}(s|m)^2+3 (m-2) q^2+\tau_j\right)
\end{gather*}
and
\[r_j = \frac{\dot{S}_{j}^{3}}{S^{3}_{j}}+\frac{S^{2}_{j}}{S^{3}_{j}} =\frac{9-9 m \text{cn}( s|m) \text{dn}(s|m) \text{sn}(s|m)}{9 \text{dn}( s|m)^2+3 (m-2) +\tau_j},\]
The quadratures can be carried out in terms of elliptic integrals of the third kind and we obtain
\begin{gather}\label{3.19}
\rho_j=\sqrt{2} \sqrt{9 \text{dn}( s|m)^2+3 (m-2) +\tau_j }\cdot e^{9\frac{ \Pi \left(\frac{9 m }{3 (m+1) +\tau_j
};\text{am}( s|m)|m\right)}{ \left(3 (m+1) +\tau_j \right)}},
\end{gather}
where
\[
\Pi(\zeta;\phi|m)=\int_0^{\phi}\frac{d\theta}{(1-\zeta\sin^2(\theta))\sqrt{1-m\sin^2(\theta)}}
\]
is the incomplete integral of the third kind and $\mathrm{am}(u|m)$ is the amplitude of the Jacobi elliptic functions. Note that on the right hand side of (\ref{3.19}) we have a smooth real branch of a multi-valued analytic function. The f\/irst column vector of~$S^{-1}$ is the transpose of
\[
\left(\frac{1}{(\tau_0-\tau_1)(\tau_0-\tau_2)},
\frac{1}{(\tau_1-\tau_0)(\tau_1-\tau_2)},\frac{1}{(\tau_2-\tau_0)(\tau_2-\tau_1)}\right)
\]
and the coef\/f\/icients $\tilde{M}_i^j$ of $M(0)^{-1}$ are
\begin{gather*}
\tilde{M}_j^1  = \frac{-9 (m+1)^2+3 (m+1) \tau_1+2 \tau_j^2}{2 \sqrt{6 (m+1) +2 \tau_j}},\\
\tilde{M}_j^2  = \frac{81}{\sqrt{6 (m+1)+2 \tau_j}},\qquad
\tilde{M}_j^3  = \frac{9 \sqrt{3 (m+1)+\tau_j}}{\sqrt{2}}.
\end{gather*}
Then, the homogeneous components of a congruence curve with projective curvature $\tilde{k}_{m}$ are
\[
\tilde{x}^j(s)=\sum_{k=0}^{2}\tilde{M}^j_k\frac{\rho_k(s)}{\prod_{h\neq k}(\tau_k-\tau_h)},\qquad j=0,1,2,
\]
and the evolution of the congruence curve is given by
\[x^j(s,t)=\sum_{k=0}^{2}\mathrm{Exp}(t\xi)^j_k \tilde{x}^k(s-(1+m(m-1)) t),\qquad j=0,1,2.\]

\begin{figure}[t]
\centering
\includegraphics[width=4cm]{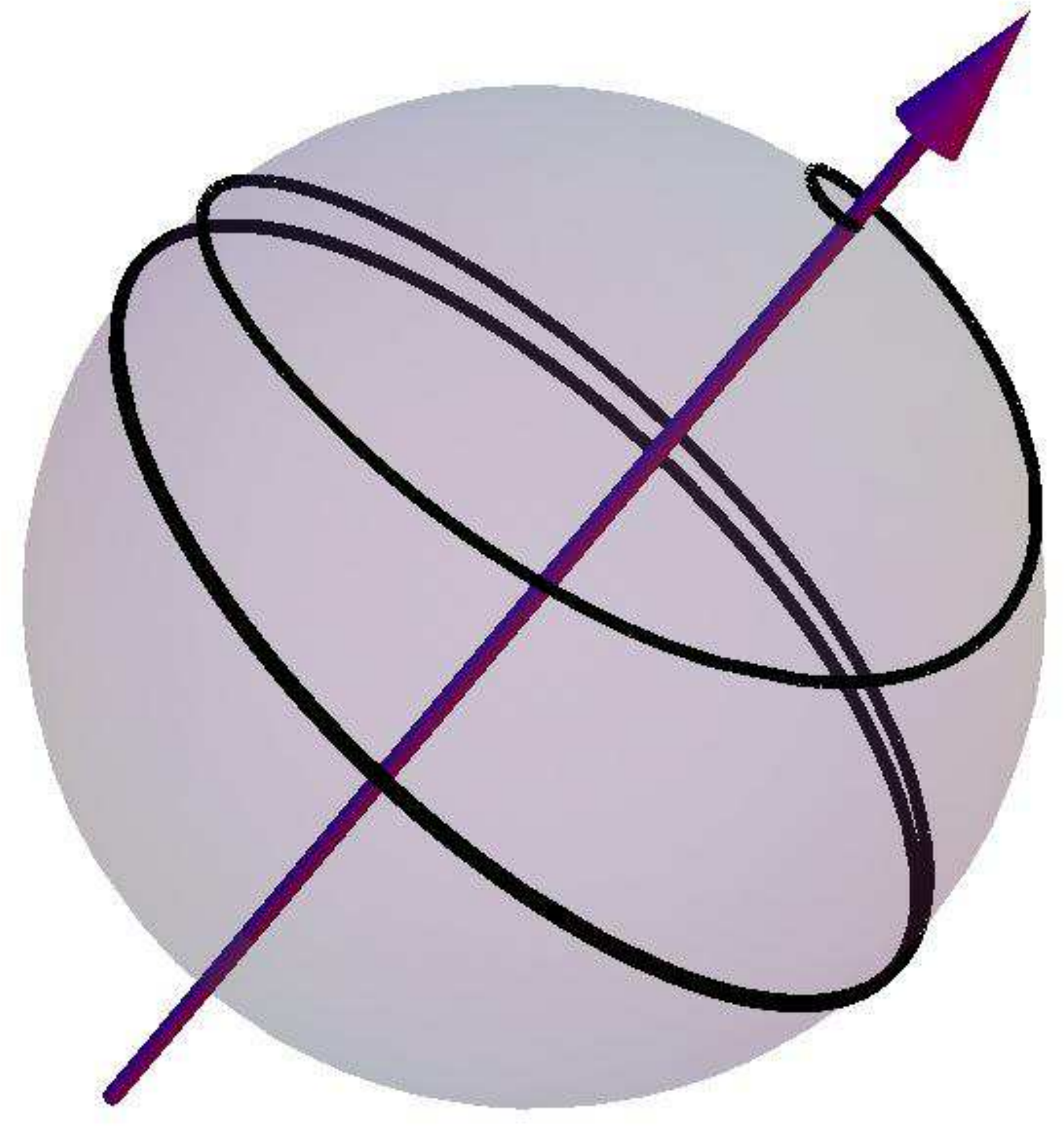}\qquad \qquad
\includegraphics[width=4cm]{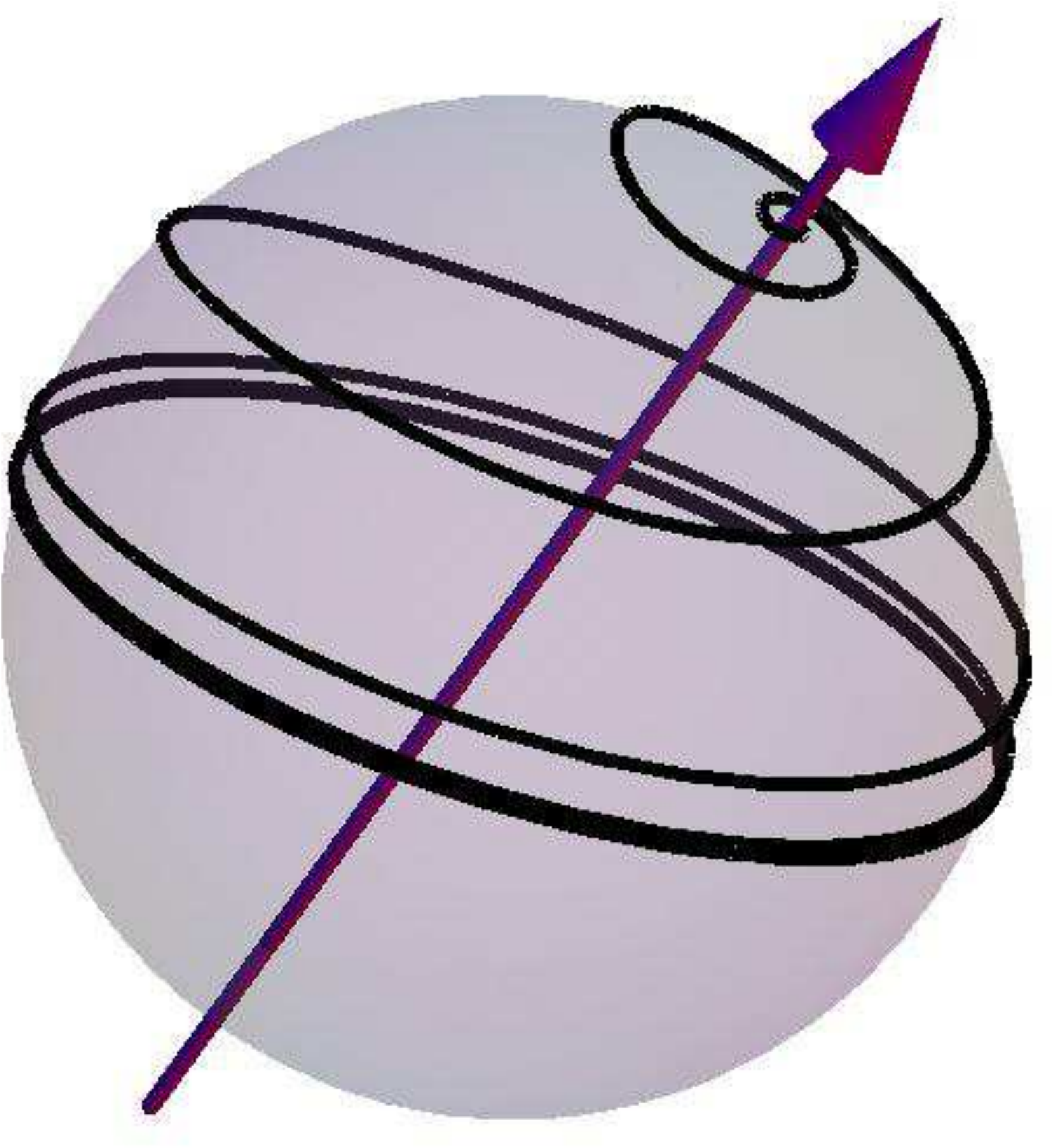}
\caption{Congruence curves with parameter $m=0.3$ and spectra $\sigma_1$ and $\sigma_2$ respectively.}\label{FIG6}
\end{figure}

\begin{remark} These curves have a spiral behavior and a gnomonic growth (i.e.\ made of successive self-congruent parts). Fig.~\ref{FIG6} reproduces the spherical lifts of the congruence curves with parameter $m=0.3$ and spectra
 \begin{gather*}
 \sigma_1  =(-6.00233 + 6.06928i,-6.00233 - 6.06928i,12.0047),\\
  \sigma_2 =(-9.94554,-9.94554,19.8911),
  \end{gather*}
respectively. Fig.~\ref{FIG7} reproduces the spherical lifts of the trajectories $\gamma(-,t)$ of the motion of a~congruence curves with parameter $m=0.7$, spectrum
\[
\sigma=(-5.45852 + 6.40263 i,10.917,-5.45852 - 6.40263 i)
\]
and $t=0,0.5,1,1.5$ respectively.
 \end{remark}

\begin{figure}[t]
\centering
\includegraphics[width=4cm]{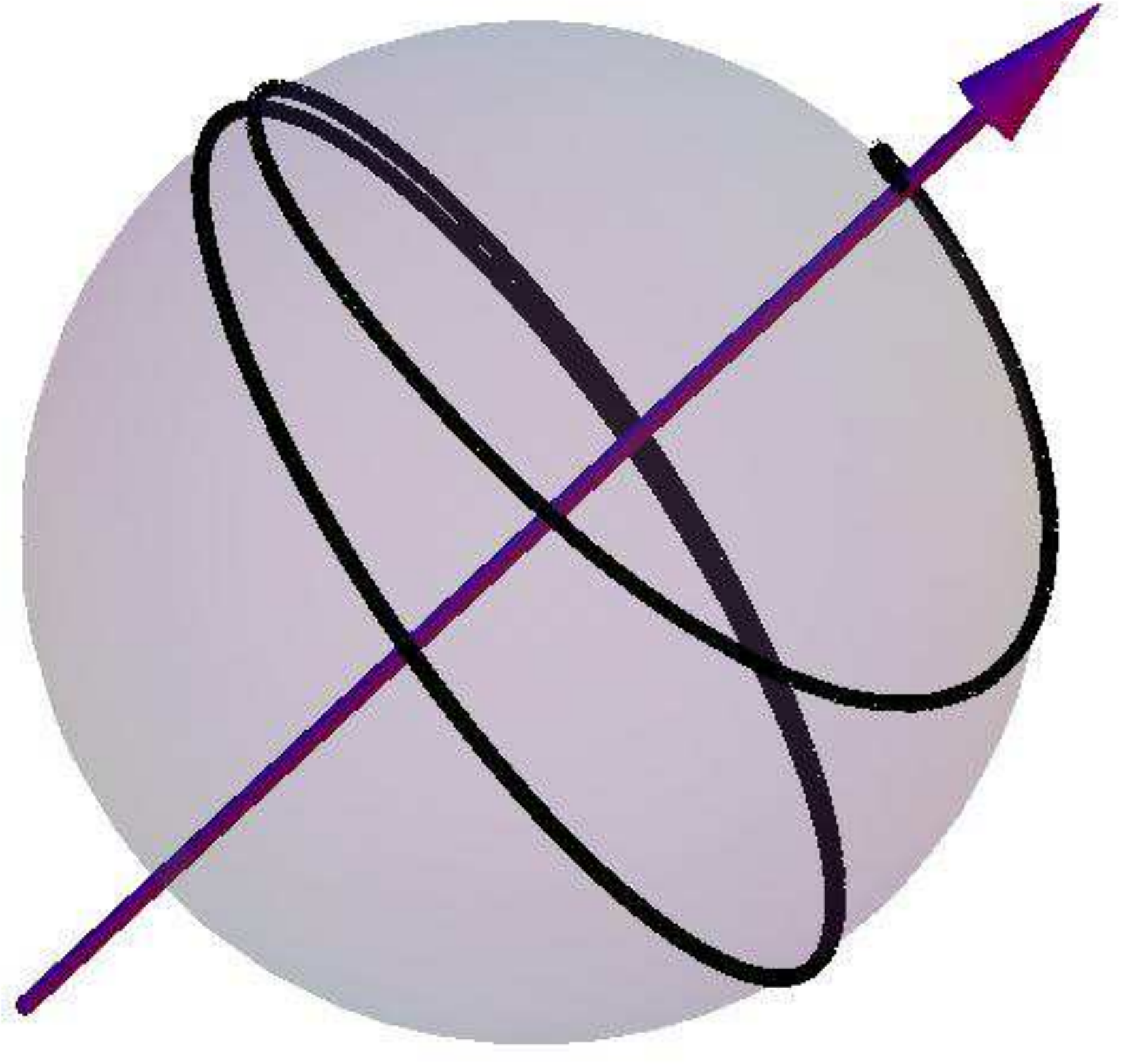}\qquad\qquad
\includegraphics[width=4cm]{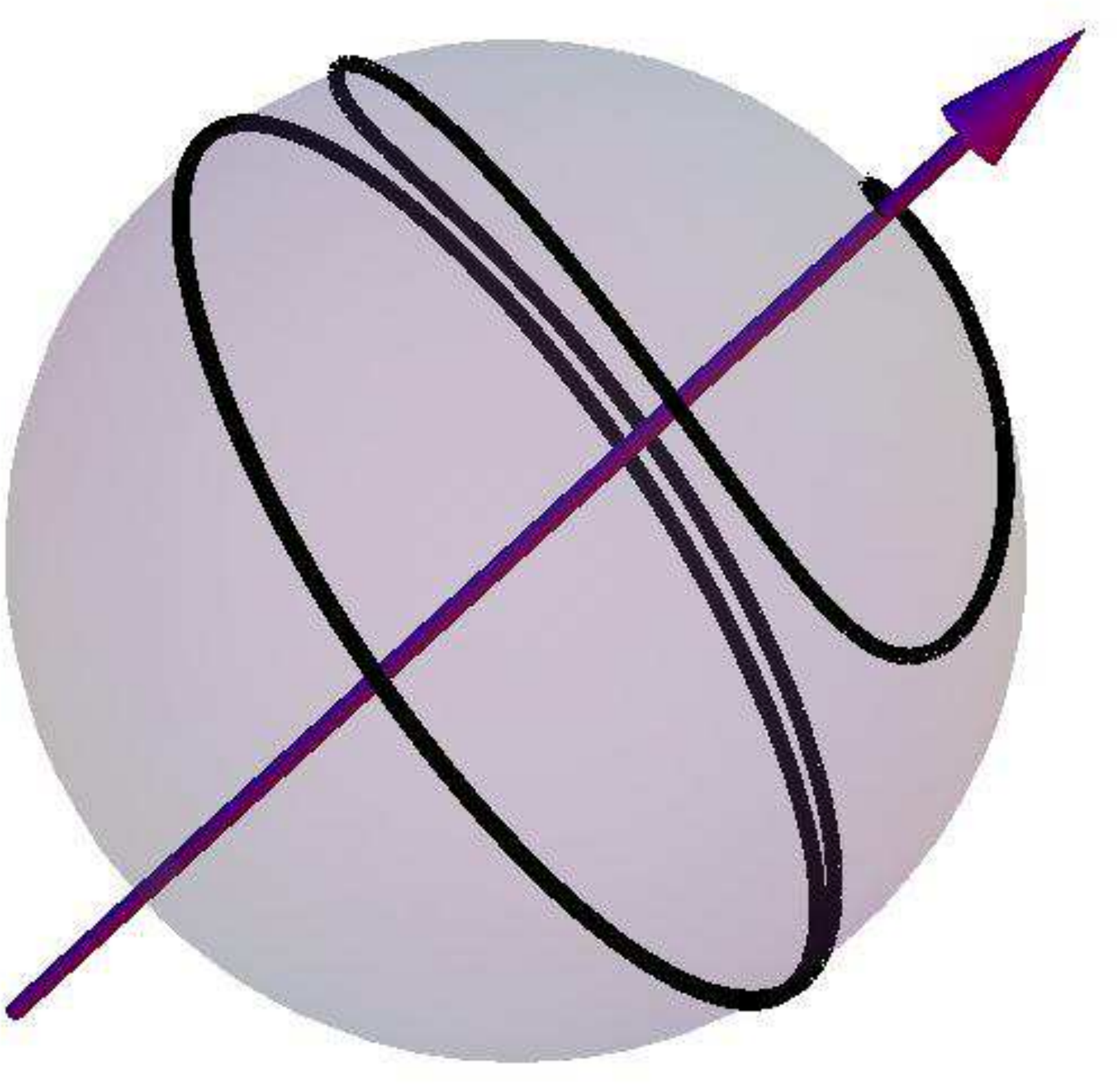}

\includegraphics[width=4cm]{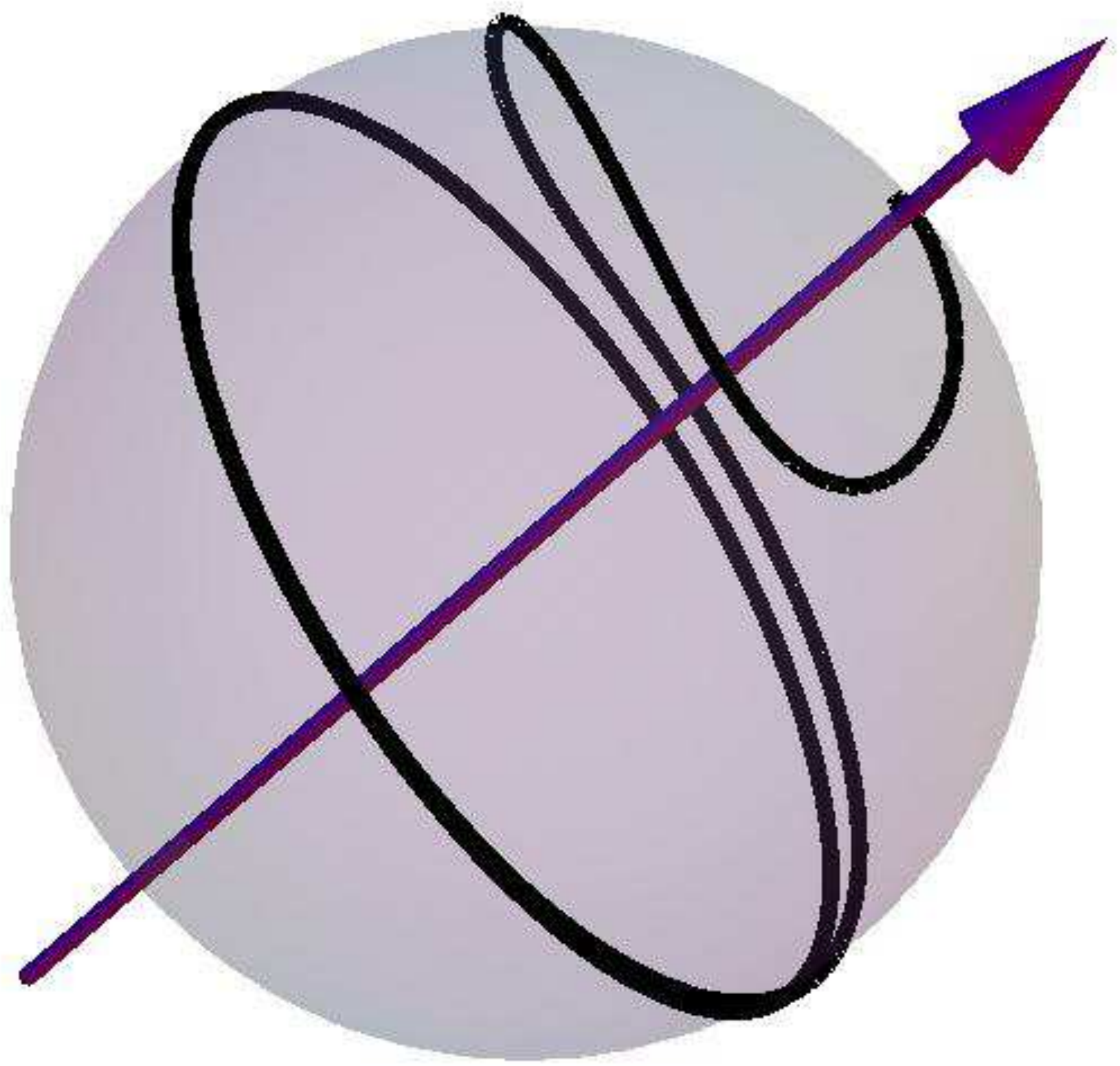}\qquad\qquad
\includegraphics[width=4cm]{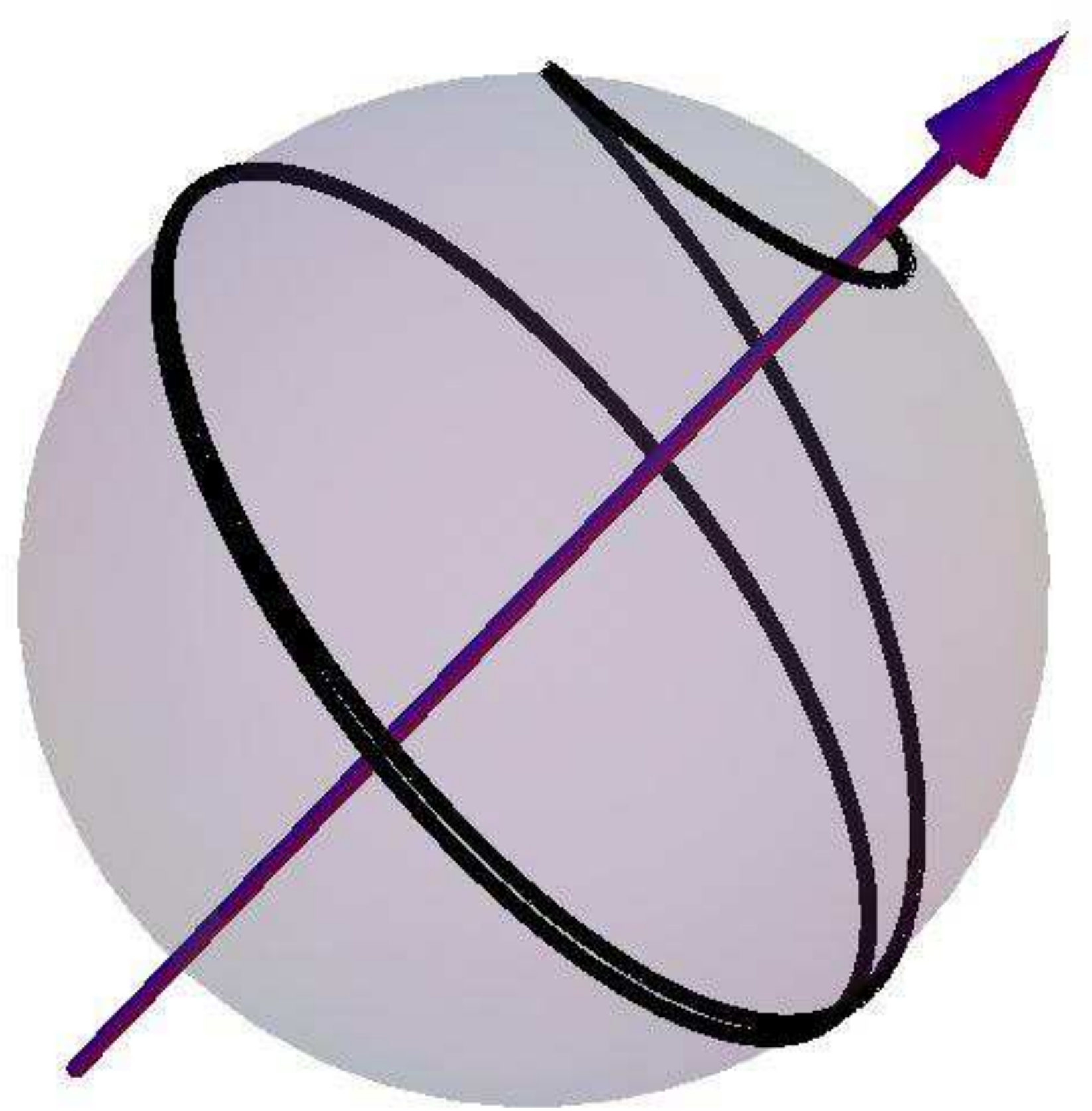}
\caption{Trajectories $\gamma(-,t)$ of the congruence curve with parameter $m=0.3$, spectrum $\sigma$ for $t=0,0.5,1$ and $1.5$ respectively.}\label{FIG7}
\end{figure}
 \end{example}

\appendix

\section{Appendix}\label{a}

\subsection[Code to compute the equations of the Kaup-Kupershmidt hierarchy]{Code to compute the equations of the Kaup--Kupershmidt hierarchy}
\label{a1}

\noindent\(\text{J1}[\text{h$\_$},\text{v$\_$}]\text{:=}D[h,\{s,3\}];\\
\text{J2}[\text{h$\_$},\text{v$\_$}]\text{:=}2*D[v*\text{Integrate}[h,s],\{s,2\}];\\
\text{J3}[\text{h$\_$},\text{v$\_$}]\text{:=}8*v{}^{\wedge}2*\text{Integrate}[h,s]+3(v*D[h,s]+D[v*h,s]);\\
\text{J4}[\text{h$\_$},\text{v$\_$}]\text{:=}2*\text{Integrate}[v*D[h,\{s,2\}]+4*v{}^{\wedge}2*h,s];\\
\mathcal{J}[\text{h$\_$},\text{v$\_$}]\text{:=}\text{Expand}[\text{FullSimplify}[\text{J1}[h,v]+\text{J2}[h,v]+\text{J3}[h,v]+\text{J4}[h,v]]];\\
\mathcal{D}[\text{h$\_$},\text{v$\_$}]\text{:=}D[h,\{s,3\}]+v*D[h,s]+D[v*h,s];\\
H[0][\text{v$\_$}]\text{:=}1;H[1][\text{v$\_$}]\text{:=}D[v,\{s,2\}]+4*v{}^{\wedge}2;H[\text{n$\_$}][\text{v$\_$}]\text{:=}
\mathcal{J}[\mathcal{D}[H[n-2][v],v],v];\\
\mathfrak{h}[\text{n$\_$}]\text{:=}\text{Expand}[H[n][u[s,t]]];\\
\mathfrak{q}[\text{n$\_$}]\text{:=}\text{Expand}[\text{Integrate}[\mathcal{D}[\mathfrak{h}[n],u[s,t]],s]];\\
\mathfrak{p}[\text{n$\_$}]\text{:=}\text{Expand}[\text{Integrate}[H[n][\epsilon *u[s,t]]*u[s,t],\{\epsilon ,0,1]]\\
\text{KK}[\text{n$\_$}]\text{:=}\text{Expand}[D[u[s,t],t]+\mathcal{D}[H[n][u[s,t]],u[s,t]]];\)

\subsection{Code to solve numerically the projective Frenet system}\label{a2}

  \emph{Step I}: def\/ine the speed, the curvature and domain of def\/inition

\noindent\(m\text{:=}0.8;\quad t\text{:=}0;\quad v[\text{s$\_$}]\text{:=}1;\quad a\text{:=}-20;\quad b\text{:=}20;\\
k[\text{s$\_$}]\text{:=}\frac{3 m^2 (1+2 \text{Cosh}[m *(s-m{}^{\wedge}4*t)])}{2 (2+\text{Cosh}[m *(s-m{}^{\wedge}4*t)])^2};\)
\medskip

\noindent \emph{Step II}: the routine to integrate the linear system

\noindent\(
\text{sol}[1]\text{:=}\text{NDSolve}\left[\left\{x'[t]==y[t],x[0]==1,y'[t]==-k[t]*x[t]+z[t],y[0]==0,\right.\right.\\
z'[t]==x[t]-k[t]*y[t],z[0]==0\},\{x,y,z\},\{t,a,b\}];\\
\text{sol}[2]\text{:=}\text{NDSolve}\left[\left\{x'[t]==y[t],x[0]==0,y'[t]==-k[t]*x[t]+z[t],y[0]==1,\right.\right.\\
z'[t]==x[t]-k[t]*y[t],z[0]==0\},\{x,y,z\},\{t,a,b\}];\\
\text{sol}[3]\text{:=}\text{NDSolve}\left[\left\{x'[t]==y[t],x[0]==0,y'[t]==-k[t]*x[t]+z[t],y[0]==0,\right.\right.\\
z'[t]==x[t]-k[t]*y[t],z[0]==1\},\{x,y,z\},\{t,a,b\}];\\
S[1][\text{t$\_$}]\text{:=}\text{Evaluate}[\{x[t],y[t],z[t]\}\text{/.}\text{sol}[1]];\\
S[2][\text{t$\_$}]\text{:=}\text{Evaluate}[\{x[t],y[t],z[t]\}\text{/.}\text{sol}[2]];\\
S[3][\text{t$\_$}]\text{:=}\text{Evaluate}[\{x[t],y[t],z[t]\}\text{/.}\text{sol}[3]];\\
\Gamma [\text{t$\_$}]\text{:=}\{S[1][t][[1]][[1]],S[2][t][[1]][[1]],S[3][t][[1]][[1]]\};\\
\gamma [\text{t$\_$}]\text{:=}\frac{1}{\sqrt{\Gamma [t].\Gamma [t]}}\Gamma [t];\)

\subsection*{Acknowledgements}

The work was partially supported by MIUR project: \textit{Metriche
riemanniane e variet\`a dif\-fe\-ren\-zia\-bi\-li}; by the GNSAGA of INDAM and by TTPU University in Tashkent.
The author would like to thank the referees and G.~Mar\'i Bef\/fa  for their useful comments and suggestions.

\pdfbookmark[1]{References}{ref}
\LastPageEnding

\end{document}